\documentclass[11pt,reqno]{amsart}
\usepackage{amsmath,amsthm,amssymb,amsfonts}
\usepackage{graphicx,xcolor,color,dsfont}
\usepackage{comment}
\usepackage{verbatim,fancyhdr,pictexwd,dcpic}
\usepackage{mathrsfs,latexsym,mathtools}
\usepackage{bm,enumerate}
\usepackage[utf8]{inputenc}
\usepackage[colorlinks=true,
    citecolor=black,
    linkcolor=black,
    urlcolor=blue,backref=page]{hyperref}

\newcommand{\strutstretchdef}{\newcommand{\strutstretch}{\vphantom{\raisebox{1pt}{$\big($}\raisebox{-1pt}{$\big($}}}}
\theoremstyle{plain}
\newtheorem{theorem}{Theorem}[section]
\newtheorem{c-theorem}{Construction theorem}[section]

\newtheorem{lemma}[theorem]{Lemma}
\newtheorem{proposition}[theorem]{Proposition}
\newtheorem{corollary}[theorem]{Corollary}

\theoremstyle{definition}
\newtheorem{definition}[theorem]{Definition}

\newtheorem{example}[theorem]{Example}
\newtheorem{question}[theorem]{Question}

\newtheorem{remark}[theorem]{Remark}

\numberwithin{equation}{section}

\newlength{\struh}
\newlength{\textminustop}
 \strutstretchdef
\makeatletter
\@namedef{subjclassname@2020}{%
\textup{2020} Mathematics Subject Classification}
\makeatother

\newcommand{\ncom}{\newcommand}
\ncom{\eop}{\hfill{{\rule{2.5mm}{2.5mm}}}}
\ncom{\eof}{\hfill{{\rule{1.5mm}{1.5mm}}}}
\ncom{\nno}{\nonumber}
\ncom{\rar}{\rightarrow}
\ncom{\Rar}{\Rightarrow}
\ncom{\noin}{\noindent}
\ncom{\sz}{\scriptsize}
\ncom{\rf}{\ref}
\ncom{\sgm}{\sigma}
\ncom{\Sgm}{\Sigma}
\ncom{\dt}{\delta}
\ncom{\Dt}{Delta}
\ncom{\lmd}{\lambda}
\ncom{\Lmd}{\Lambda}
\ncom{\eps}{\epsilon}
\ncom{\pcc}{\stackrel{P}{>}}
\ncom{\dist}{{\rm\,dist}}
\ncom{\im}{{\rm Im\,}}
\ncom{\sgn}{{\rm sgn\,}}
\ncom{\hone}{\mbox{\hspace{1em}}}
\ncom{\htwo}{\mbox{\hspace{2em}}}
\ncom{\hthree}{\mbox{\hspace{3em}}}
\ncom{\hfour}{\mbox{\hspace{4em}}}
\ncom{\hsev}{\mbox{\hspace{7em}}}
\ncom{\vone}{\vskip 2ex}
\ncom{\vtwo}{\vskip 4ex}
\ncom{\vonee}{\vskip 1.5ex}
\ncom{\vthree}{\vskip 6ex}
\ncom{\vfour}{\vspace*{8ex}}
\ncom{\norm}{\|\;\;\|}
\ncom{\integ}[4]{\int_{#1}^{#2}\,{#3}\,d{#4}}
 \ncom{\inp}[2]{\langle{#1},\,{#2} \rangle}
\ncom{\Inp}[2]{\Langle{#1},\,{#2} \Langle}
\ncom{\vspan}[1]{{{\rm\,span}\#1 \}}}
\ncom{\dm}[1]{\displaystyle {#1}}

\begin{document}
\title[Gleason's problem and Left-inverse commuting submodules]{Gleason's problem and left-inverse commuting submodules of  Dirichlet-type spaces on  the unit bidisc}

\author[S. Bera, S. Ghara]{Santu Bera and Soumitra Ghara}

\address{Department of Mathematics\\ Indian Institute of Technology Kharagpur, 721302, India}
\email{santu@maths.iitkgp.ac.in, santuddd12@gmail.com}
\email{soumitra@maths.iitkgp.ac.in}

\thanks{The research of the first named author is supported by the Institute Postdoctoral Fellowship (ID: PDF2026176) at the Indian Institute of Technology Kharagpur.  The work of the second named author is supported by ARG-MATRICS grant by the ANRF (File No: ANRF/ARGM/2025/001583/MTR) and INSPIRE Faculty Fellowship (DST/INSPIRE/04/2021/002555).}



\begin{abstract}  
We consider harmonically weighted Dirichlet spaces $\mathcal D(\mu_1,\mu_2)$ on the unit bidisc $\mathbb D^2$ induced by two finite nonnegative Borel measures $\mu_1$ and $\mu_2$ supported on the unit circle $\mathbb T$. These spaces provide a unified framework for studying  the Hardy and Dirichlet spaces on the bidisc. Moreover, these spaces serve as model spaces for a class of cyclic analytic toral $2$-isometries. Our contributions in this paper are twofold.

Our first contribution is a solution to Gleason's problem for both $\mathcal D(\mu_1,\mu_2)$ and its multiplier algebra. As applications, we prove that the Taylor essential spectrum of the coordinate multiplication tuple $(M_{z_1},M_{z_2})$  is $(\mathbb T\times \overline{\mathbb D})\cup(\overline{\mathbb D}\times \mathbb T)$ and  establish an index theorem for invariant subspaces which contain a nonzero multiplier.

Our second result concerns the structure of a class of submodules of $\mathcal D(\mu_1,\mu_2).$  On the Hardy space $H^2(\mathbb D^2)$, the  multiplication tuple $(M_{z_1}, M_{z_2})$ is doubly commuting, whereas on $\mathcal D(\mu_1,\mu_2)$ it is left-inverse commuting. Motivated by Mandrekar's characterization of doubly commuting submodules of $H^2(\mathbb D^2)$, we establish an analogous characterization of left-inverse commuting submodules of $\mathcal D(\mu_1,\mu_2)$. In the process, we also show that every nonzero left-inverse commuting submodule $\mathcal S$ satisfies
$$\dim\bigl(\mathcal S\ominus(z_1\mathcal S+z_2\mathcal S)\bigr)=1.$$
Furthermore, in the case of the Dirichlet space on \(\mathbb D^2\), we show that the nonzero left-inverse commuting submodules are precisely those generated by monomials.
\end{abstract}

\maketitle
\tableofcontents
\section{Introduction}
\subsection{Left-inverse commuting submodules of \texorpdfstring{$\mathcal D(\mu_1,\mu_2)$}{}:}
In operator theory on analytic function spaces, a central problem is to describe all closed invariant subspaces of a given operator
Beurling in his seminal paper \cite{B1949} proved that every closed invariant subspace $\mathcal S$ of the Hardy space $H^2(\mathbb D)$ of the unit disc $\mathbb D$ satisfies $\dim (\mathcal S\ominus z\mathcal S)=1$ and  is of the form $qH^2(\mathbb D)$ for some inner function $q$ on $\mathbb D.$ On the other hand, Richter characterised the invariant subspaces of the Dirichlet space of the unit disc. Recall that the Dirichlet space $\mathcal D(\mathbb D)$ on the unit disc is defined by
\begin{equation*}
\mathcal{D}(\mathbb D):=\{f\in \mathcal O(\mathbb D):\int_{\mathbb D}|f'(z)|^2 dA(z)<\infty\},    
\end{equation*}
with the norm 
$$\|f\|_{\mathcal D(\mathbb D)}^2=\|f\|_{H^2(\mathbb D)}^2+\int_{\mathbb D}|f'(z)|^2 dA(z),$$
where $dA$ denotes the \textit{normalized Lebesgue area measure} on the unit disc, and   
for a domain $\Omega\subset \mathbb C^d$ $(d\geqslant 1)$, $\mathcal O(\Omega)$ denotes the space of holomorphic functions on $\Omega$. To characterise  invariant subspaces of $\mathcal D(\mathbb D)$, Richter \cite{R1988} introduced a notion of {\it Dirichlet-type spaces} $\mathcal D(\mu)$ on unit disc as follows: for a finite nonnegative Borel measure $\mu$ supported on the unit circle  $\mathbb T,$

\begin{equation*}
\mathcal D(\mu):=\{f\in H^2(\mathbb D):  D_{\mu}(f):=\int_{\mathbb D} |f'(z)|^2 P_\mu(z) dA(z)<\infty\},
\end{equation*}
where 
$$P_{\mu}(z)=\int_{\mathbb T}\frac{1-|z|^2}{|z-\zeta|^2}d\mu(\zeta).$$ A routine verification shows that $ \mathcal D(\mu)$ is a reproducing kernel Hilbert space with respect to the norm 
$$\|f\|_{\mathcal D(\mu)}^2:=\|f\|_{H^2(\mathbb D)}^2+ D_{\mu}(f).$$
Note that if $m$ denotes the \textit{normalised Lebesgue measure} on $\mathbb T,$  $\mathcal D(m)$ is the Dirichlet space of $\mathbb D$. It is observed that these spaces serve as model spaces for cyclic analytic $2$-isometries \cite[Theorem~5.1]{R1991}. In fact, any closed invariant subspace $\mathcal S$ of $\mathcal D(\mu)$  satisfies $\dim (\mathcal S\ominus z\mathcal S)=1$ and  is of the form $\varphi\mathcal D(\mu_\varphi)$ for some unit vector $\varphi \in  \mathcal S\ominus z\mathcal S$, where $d\mu_\varphi(\zeta)= |\varphi^*(\zeta)|^2d\mu(\zeta)$ (see \cite[Theorem~7.1]{R1988} and \cite[Theorem~3.2]{RS1992}) (we refer the reader  to \cite{EKMR2014} for a  comprehensive study of harmonically weighted Dirichlet spaces on $\mathbb D$). Any subspace of the similar form is often referred as a {\it Beurling type invariant subspace}.

However, the situation changes substantially when one moves from one variable to several complex variables. To illustrate this phenomenon, let $H^2(\mathbb D^2)$ denote the Hardy space on the bidisc $\mathbb D^2$, defined by
$$ H^2(\mathbb D^2):=\left\{f\in\mathcal O(\mathbb D^2):\|f\|_{H^2(\mathbb D^2)}^2:=\sup_{0\leq r<1} \int_0^{2\pi}\int_0^{2\pi}|f(re^{i\theta},re^{i\phi})|^2\frac{d\theta}{2\pi} \frac{d\phi}{2\pi}<\infty\right\}.$$
In sharp contrast with the one-variable case, Rudin \cite{Ru1969} showed that not every closed invariant subspace of $H^2(\mathbb D^2)$ admits a Beurling-type representation. More strikingly, he constructed a closed invariant subspace that contains no nonzero bounded holomorphic function. This reveals the considerably richer structure of invariant subspaces in several variables and naturally leads to the problem of characterizing those invariant subspaces that retain a Beurling-type structure.

A fundamental contribution to this problem was made by Mandrekar \cite[Theorem~2]{M1988} (see also \cite{GM1988}), who established that a nonzero closed invariant subspace $\mathcal S$ of $H^2(\mathbb D^2)$ can be represented as $\mathcal S=qH^2(\mathbb D^2)$ for some inner function $q$ on $\mathbb D^2$ if and only if the pair of restricted coordinate multiplication operators $(M_{z_1}|_{\mathcal S},M_{z_2}|_{\mathcal S})$ is doubly commuting. Recall that a commuting pair $(T_1,T_2)$ of operators on a Hilbert space is said to be {\it doubly commuting} if $T_1^*T_2=T_2T_1^*$. Also, a bounded holomorphic function $q$ on $\mathbb D^2$ is called {\it inner} if its radial boundary values satisfy
$$|q^*(\zeta_1,\zeta_2)|=1 \quad\text{for almost every }(\zeta_1,\zeta_2)\in\mathbb T^2.$$
Thus, Mandrekar's theorem establishes a remarkable connection between the function-theoretic structure of invariant subspaces and the operator-theoretic property of doubly commutativity.

 Motivated by Richter's notion of $\mathcal D(\mu)$-spaces, a notion of Dirichlet-type spaces $\mathcal D(\mu_1,\mu_2)$ on $\mathbb D^2$ was introduced recently in \cite{BCG2025}. These spaces have attracted considerable attention very recently, see \cite{B2026,BS2026,BGV2025,BGV2026,NZ2026}.
 For any $f\in \mathcal O(\mathbb D^2)$ and two finite  nonnegative Borel measures $\mu_1, \mu_2$, define $D_{\mu_1,\mu_2}(f)$ by 
$$D_{\mu_1, \mu_2}(f)=D^{(1)}_{\mu_1, \mu_2}(f)+D^{(2)}_{\mu_1, \mu_2}(f),$$
where 
 \begin{align*}
D^{(1)}_{\mu_1, \mu_2}(f) &:= \lim_{r\to 1^-} \int_{0}^{2\pi} \int_{\mathbb D}|\partial_1 f(z_1, re^{i \theta})|^2 P_{\mu_1}(z_1) \,dA(z_1)\frac{d\theta}{2\pi}, \\
D^{(2)}_{\mu_1, \mu_2}(f)&:= \lim_{r\to 1^-} \int_{0}^{2\pi} \int_{\mathbb D}|\partial_2 f(re^{i \theta}, z_2)|^2 P_{\mu_2}(z_2) \,dA(z_2)\frac{d\theta}{2\pi}. 
\end{align*} 
The \emph{Dirichlet-type space} $\mathcal D(\mu_1,\mu_2)$ is defined as the space of functions $f\in H^2(\mathbb D^2)$ for which 
$$\|f\|^2_{\mathcal D(\mu_1, \mu_2)}:=\|f\|^2_{H^2(\mathbb D^2)}+D_{\mu_1,\mu_2}(f)<\infty.$$
It turns out that $\mathcal D(\mu_1,\mu_2)$ is a reproducing kernel Hilbert space 
on $\mathbb D^2$ and the set of polynomials is dense in it. Moreover,  
$\mathcal D(\mu_1,\mu_2)$ serves as a model space for a class of cyclic analytic toral $2$-isometries (see \cite[Theorem~2.4]{BCG2025}). Recall that a pair $T=(T_1,T_2)$ of commuting bounded linear operators on a complex separable Hilbert space $\mathcal H$ is said to be a {\it toral $2$-isometry} if
\begin{equation}
\label{C1}
I - T_i^*T_i - T_j^*T_j+T_i^*T^*_jT_iT_j=0, \quad i,j = 1,2
\end{equation} 
(see \cite{A1990, AS1999, BCG2025, R1991}). Note that when $\mu_1=\mu_2=0$, then  $\mathcal D(\mu_1,\mu_2)$ coincides with the Hardy space $H^2(\mathbb D^2)$. Also, the space $\mathcal D(m,m)$ is the Dirichlet space on $\mathbb D^2$ defined by
\begin{align*}
\mathcal D(\mathbb D^2):=\Big\{f =\sum_{k,\ell=0}^\infty \hat{f}(k,\ell) z_1^k z_2^\ell\in \mathcal O(\mathbb D^2): \sum_{k,\ell=0}^\infty (k+\ell)|\hat{f}(k,\ell)|^2  < \infty\Big\},
\end{align*}
 with the norm 
\begin{align*}
\|f\|^2_{\mathcal D(\mathbb D^2)}:=\sum_{k,\ell=0}^\infty|\hat{f}(k,\ell)|^2 (k+\ell+1),\quad f\in \mathcal D(\mathbb D^2).
\end{align*}

Following \cite{BGV2025}, a commuting pair $T=(T_1,T_2)$ of left invertible operators is called {\it left-inverse commuting} if $L_iT_j=T_jL_i$ for $1\leqslant i\neq j\leqslant 2,$ where $L_i$ denotes the left inverse of $T_i$ given by $$L_i=(T_i^*T_i)^{-1}T_i^*,\quad i=1,2.$$ Clearly, if  $T=(T_1, T_2)$ is a commuting pair of isometries, then $T$ is left-inverse commuting if and only if $T$ is doubly commuting. Furthermore, it is observed in \cite[Example~5.3]{BGV2025} that 
\begin{equation}\label{eqnlic}
\text{the multiplication tuple 
$M_z=(M_{z_1},M_{z_2})$ on $\mathcal D(\mu_1,\mu_2)$ is left-inverse commuting.} 
\end{equation}
A closed subspace  $\mathcal S$ of $\mathcal D(\mu_1,\mu_2)$ is called a {\it submodule} of $\mathcal D(\mu_1,\mu_2)$ if $\mathcal S$ is invariant under both the multiplication operators $M_{z_1}$ and $M_{z_2}$.  A submodule $\mathcal S$ of $\mathcal D(\mu_1,\mu_2)$ is called left-inverse commuting if the pair $(M_{z_1}|_{\mathcal S},M_{z_2}|_{\mathcal S})$ is left-inverse commuting. In view of \eqref{eqnlic} and motivated by Mandrekar's characterisation of doubly-commuting submodules of $H^2(\mathbb D^2)$, we ask the following question:
\begin{question}\label{ques}
        What are all left-inverse commuting submodules of $\mathcal D(\mu_1,\mu_2)?$
\end{question}

Recall that a function $\psi\in \mathcal D(\mu)$ is called  extremal in $\mathcal D(\mu)$ if $\|\psi\|_{\mathcal D(\mu)}=1$ and $\langle \psi, z^n\psi \rangle_{\mathcal D(\mu)}=0$ for all $n \geqslant 1.$ These functions play an important role in the classification of invariant subspaces of $\mathcal D(\mu)$ spaces. In fact, every nonzero invariant subspace of $\mathcal D(\mu)$ is generated by an extremal function, 
(see \cite[Theorem~3.1]{RS1992}, \cite[Theorem~IV.4.9]{A1993}).
Note that $\psi$ is extremal in $H^2(\mathbb D)$ if and only if it is an inner function. Motivated by this, 
we define the following notion of extremal function in $\mathcal D(\mu_1,\mu_2)$:
\begin{definition}
     A  function $\varphi\in \mathcal D(\mu_1,\mu_2)$ is said to be extremal in $\mathcal D(\mu_1,\mu_2)$  if   $\|\varphi\|_{\mathcal D(\mu_1,\mu_2)} = 1$ and 
    \begin{equation}\label{extremal}
        \begin{split}
            \langle z_1^k\varphi, z_1^pz_2^{q+1}\varphi\rangle_{\mathcal D(\mu_1,\mu_2)}&=0,\quad k,p,q\geqslant 0,\\
            \langle z_2^\ell\varphi, z_1^{p+1}z_2^q\varphi\rangle_{\mathcal D(\mu_1,\mu_2)}&=0,\quad  \ell,p,q\geqslant 0.
        \end{split}
    \end{equation}     
         \end{definition}
It is not difficult to see that the only extremal functions in $H^2(\mathbb D^2)$ are inner functions on $\mathbb D^2$.  

In the following theorem, we provide a complete solution to Question \ref{ques}.

\begin{theorem}\label{inv-sub-d-mu}
Let $\mu_1, \mu_2$ be two finite nonnegative Borel measures on $\mathbb T$. The following statements hold:
\begin{enumerate}
   \item[\rm(i)] Every nonzero left-inverse commuting submodule  $\mathcal S$  of $\mathcal D(\mu_1, \mu_2)$ satisfies  $$\dim \mathcal S\ominus (z_1\mathcal S+z_2\mathcal S)=1.$$
\item[\rm(ii)] Let $\mathcal S$ be a nonzero submodule of  $\mathcal D(\mu_1, \mu_2)$. 
Then $\mathcal S$ is left-inverse commuting if and only if $\mathcal S$ is of the form $\varphi \mathcal D(\nu_1,\nu_2)$ for some extremal function $\varphi $ in  $\mathcal D(\mu_1,\mu_2)$ and finite nonnegative Borel measures $\nu_1,\nu_2$ on $\mathbb T$ of the form 
\begin{equation}\label{nu-1-2-formula}
d\nu_1(\zeta)=\sum_{\ell=0}^\infty |\varphi^{(1)*}_\ell(\zeta)|^2 d\mu_1(\zeta),\quad 
d\nu_2(\zeta)=\sum_{k=0}^\infty |\varphi^{(2)*}_k(\zeta)|^2 d\mu_2(\zeta),
\end{equation}
where $\varphi^{(1)}_\ell, \varphi^{(2)}_k\in H^2(\mathbb D)$, $k,\ell\in \mathbb Z_+$, are given by 
\begin{equation*}
    \varphi(z_1,z_2)=\sum_{k=0}^\infty \varphi^{(2)}_k(z_2) z_1^k=\sum_{\ell=0}^\infty \varphi^{(1)}_\ell(z_1) z_2^\ell, \quad (z_1,z_2)\in \mathbb D^2.
\end{equation*}
\end{enumerate}
\end{theorem}

In the case of $\mu_1=\mu_2=0,$ that is $\mathcal D(\mu_1,\mu_2)=H^2(\mathbb D^2),$ using the fact that extremal functions in $H^2(\mathbb D^2)$ are inner functions, the above theorem recovers Mandrekar's characterisation of doubly commuting submodules of $H^2(\mathbb D^2)$ (see \cite[Theorem~2]{M1988}). In the case of $\mathcal D(\mathbb D^2)$, we show that the extremal functions are constant multiples of monomials.
\begin{theorem}\label{extremals-are-monomials}
A function $\varphi\in\mathcal D(\mathbb D^2)$ is extremal in $\mathcal D(\mathbb D^2)$ if and only if there exist  $k,\ell\in \mathbb Z_+$ and $\theta\in\mathbb R$ such that 
\begin{equation*}
\varphi(z_1,z_2)=\frac{e^{i\theta}}{\sqrt{1+k+\ell}}z_1^kz_2^\ell.
\end{equation*}
\end{theorem}
The preceding result fails if one of $\mu_1$ or $\mu_2$ is not equal to the Lebesgue measure (see Remark~\ref{non-monomial-example}).
As an application of Theorems~\ref{inv-sub-d-mu}~\&~\ref{extremals-are-monomials} we establish the following classification of left-inverse commuting submodules of $\mathcal D(\mathbb D^2)$:
\begin{corollary}\label{inv-d-m-m}
Let   $\mathcal S\neq \{0\}$ be a submodule of $\mathcal D(\mathbb D^2)$. Then $\mathcal S$ is left-inverse commuting if and only if  $\mathcal S=[z_1^kz_2^\ell],$ the cyclic subspace generated by the monomial $z_1^kz_2^\ell,$ where $k,\ell\in \mathbb Z_+.$ 
In this case, $\mathcal S$ is of the form
$$\mathcal S=\varphi\mathcal D(cm,cm),\quad \text{where  $\varphi=\frac{z_1^kz_2^\ell}{\sqrt{k+\ell+1}}$, $c=\frac{1}{k+\ell+1}$}.$$
\end{corollary}

\subsection{Gleason's problem for \texorpdfstring{$\mathcal D(\mu_1,\mu_2)$}{}:} One of the interesting problem for spaces of holomorphic functions in several variables is the Gleason's problem. For a positive integer $d,$ let $\Omega$ be a domain in $\mathbb C^d$ and let $\mathscr B$ be a Banach space of holomorphic functions on $\Omega.$ We say that {\it Gleason's problem can be solved for $\mathscr B$} if for every $f \in \mathscr B$ and $\lambda = (\lambda_1, \ldots, \lambda_d) \in \Omega,$ there exist functions $f_1, \ldots, f_d$ in $\mathscr B$ such that
\begin{align*}
f(z)=f(\lambda) + \sum_{j=1}^d (z_j-\lambda_j)f_j(z), \quad z=(z_1, \ldots, z_d) \in \Omega.
\end{align*}
Gleason's problem is known to be solvable for many classical function spaces, for instance, the Hardy space $H^2(\mathbb D^d)$ of the unit polydisc, the Hardy and Bergman spaces of the unit ball $\mathbb B^d$ (see \cite[Corollary~4.6]{GRS2005}, \cite{Z1988}). The same is true for their multiplier algebras $H^\infty(\mathbb D^d)$ and $H^\infty(\mathbb B^d)$, respectively. Another notable example where  Gleason's problem is solvable is the Drury-Arveson space (see \cite[Section~3]{AK2002}) and its multiplier algebra (see \cite[Corollary~4.2]{GRS2005}). 

Let $\Omega$ be a domain in $\mathbb C^d.$ Let $\mathscr H$ be a reproducing kernel Hilbert space of holomorphic functions on $\Omega.$ A function $\varphi : \Omega \rightarrow \mathbb C$ is said to be a {\it multiplier} of $\mathscr H$ if $\varphi f \in \mathscr H$ for every $f \in \mathscr H.$ The space of multipliers of $\mathscr H$ is denoted by $\mathrm{Mult}(\mathscr H).$ 
The following is one of the main results of this paper.

\begin{theorem}\label{thm:gleason}
 Gleason's problem is solvable  for  $\mathcal D(\mu_1,\mu_2)$ and $\mathrm{Mult}(\mathcal D(\mu_1,\mu_2)).$ 
\end{theorem}

In \cite{GRS2005}, Gleason, Richter and Sundberg studied the Fredholm index of the multiplication tuple restricted to an invariant subspaces of $H^2(\mathbb D^d)$. In particular, \cite[Theorem~1.1]{GRS2005} shows that if an invariant subspace $\mathcal S$ of $H^2(\mathbb D^d)$ contains a  multiplier $\varphi,$ then the Taylor essential spectrum $\sigma_e(M_z|_\mathcal S)$ of the restriction tuple $M_z|_{\mathcal S}$ satisfies 
$$\sigma_e(M_z|_\mathcal S)\cap \mathbb D^d\subseteq Z(\varphi),$$
where $Z(\varphi)$ denotes the zero set of $\varphi$ inside $\mathbb D^d$ 
and for each $\lambda \in \mathbb D^d \backslash \sigma_e(M_z|_\mathcal S)$ the tuple $(M_z-\lambda)|_\mathcal S$ has Fredholm index $(-1)^d.$ In fact, for $\lambda\in \mathbb D^d\backslash Z(\varphi)$ it is shown that  
\begin{equation*}
\dim \mathcal S/\left( (z_1-\lambda_1)\mathcal S+\cdots+(z_d-\lambda_d)\mathcal S\right) =1.
\end{equation*}
Their approach relies on the solvability of Gleason's problem on both the underlying function spaces and their multiplier algebras, which yields an isomorphism between the augmented Koszul complexes associated with the multiplication tuple and its restriction. The aforementioned result is, in fact, a particular case of the more general result established in \cite[Theorem~3.4]{GRS2005}.
 Using Theorem \ref{thm:gleason},  we obtain the following result. Recall that a submodule $\mathcal S$ of $\mathcal D(\mu_1,\mu_2)$ is called \textit{multiplier invariant} if $\varphi f\in \mathcal S$ for every $f\in \mathcal S.$ 
\begin{theorem}\label{codim-thm}
Let $\mathcal S$ be a multiplier invariant submodule of $\mathcal D(\mu_1,\mu_2)$ containing a nonzero multiplier $\varphi.$ Then
\begin{align*}
\sigma_e(M_z|_\mathcal S) \cap\mathbb D^2\subseteq Z(\varphi)
\end{align*}
and for each $\lambda=(\lambda_1,\lambda_2) \in \mathbb D^2\backslash \sigma_e(M_z|_\mathcal S)$, the tuple $(M_z-\lambda)|_\mathcal S$ has Fredholm index $1.$ In fact, for all $\lambda=(\lambda_1,\lambda_2) \in \mathbb D^2\backslash Z(\varphi)$, we have 
\begin{equation*}
\dim \mathcal S/\left( (z_1-\lambda_1)\mathcal S+(z_2-\lambda_2)\mathcal S\right) =1.
\end{equation*}
\end{theorem}

The remainder of the paper is organized as follows. In Section~\ref{section-2}, we establish some preliminary results which will be used throughout this paper. In Section~\ref{section-3}, we provide proofs of Theorems \ref{thm:gleason} and \ref{codim-thm}. Section~\ref{section-4} is devoted to the proof of Theorem~\ref{inv-sub-d-mu}. In Section~\ref{section-5}, we first show that every extremal function in $\mathcal D(\mathbb D^2)$ is a constant multiple of an inner (Lemma~\ref{phi-inner}) and use this to prove Theorem~\ref{extremals-are-monomials}. Also we provide an example of a submodule of $\mathcal D(\mathbb D^2)$ that fails to be left-inverse commuting (Example~\ref{not-LIC}). Finally, we prove that $\mathcal D(\mathbb D^2)$ does not have any nonzero doubly commuting submodule (Theorem~\ref{non-d-commuting}).

\section{Left and right evaluation operators on  \texorpdfstring{$\mathcal D(\mu_1, \mu_2)$}{}}\label{section-2}
Following \cite{Y1999,Y2001}, for $\lambda\in \mathbb D,$ we define the {\it left and right evaluation operators} on $\mathcal O(\mathbb D^2)$ as
\begin{align*}
(L_\lambda f)(w):=f(\lambda,w) \text{ and } (R_\lambda f)(w):=f(w,\lambda),\quad  w\in \mathbb D,~f\in \mathcal O(\mathbb D^2).
\end{align*}
It is well-known (see \cite[pp.~382-383]{Y1999}) that if $f\in H^2(\mathbb D^2)$, then 
\begin{equation}\label{eqnhardy1}
   \text{both $L_\lambda f$ and $R_\lambda f$ belong to  $H^2(\mathbb D)$ for all $\lambda\in\mathbb D$,}
\end{equation}
 and
 \begin{equation*}
\|f\|_{H^2(\mathbb D^2)}^2=\lim_{r\to 1^-}\int_{0}^{2\pi}\|L_{re^{i\theta}} f\|_{H^2(\mathbb D)}^2\frac{d\theta}{2\pi}=\lim_{r\to 1^-}\int_{0}^{2\pi}\|R_{re^{i\theta}} f\|_{H^2(\mathbb D)}^2\frac{d\theta}{2\pi},
\end{equation*}
that is,
\begin{equation}\label{eqnhardy}
\|f\|_{H^2(\mathbb D^2)}^2=\lim_{r\to 1^-}\int_{0}^{2\pi}\|f(\cdot, re^{i\theta})\|^2_{H^2(\mathbb D)}\frac{d\theta}{2\pi}=\lim_{r\to 1^-}\int_{0}^{2\pi}\|f(re^{i\theta}, \cdot)\|^2_{H^2(\mathbb D)}\frac{d\theta}{2\pi}.
\end{equation}

The following lemma serves as a key ingredient for solving Gleason's problem for $\mathcal D(\mu_1,\mu_2)$. 
\begin{lemma}\label{slices-d-mu}
For $\lambda\in \mathbb D,$ $L_\lambda: \mathcal D(\mu_1,\mu_2) \to  \mathcal D(\mu_2)$ and $R_\lambda: \mathcal D(\mu_1,\mu_2) \to  \mathcal D(\mu_1)$ are bounded. Moreover,
\begin{align*}
\|L_\lambda f\|^2_{\mathcal D(\mu_2)}\leqslant \frac{1}{1-|\lambda|^2}\|f\|^2_{\mathcal D(\mu_1,\mu_2)},~~\|R_\lambda f\|^2_{\mathcal D(\mu_1)}\leqslant \frac{1}{1-|\lambda|^2}\|f\|^2_{\mathcal D(\mu_1,\mu_2)},~f\in\mathcal D(\mu_1,\mu_2).
\end{align*}
\end{lemma}
\begin{proof}
Let $\lambda\in\mathbb D$.  For $f\in H^2(\mathbb D^2)$, by \eqref{eqnhardy1} and \eqref{eqnhardy}, we have
\begin{align}
    \|L_{\lambda}f\|_{H^2(\mathbb D)}^2=\lim_{r\to 1^-}\int_{0}^{2\pi}|(L_{\lambda}f)(re^{i\theta})|^2\frac{d\theta}{2\pi}&=\lim_{r\to 1^-}\int_{0}^{2\pi}|f(\lambda, re^{i\theta})|^2\frac{d\theta}{2\pi} \notag\\
    &\leqslant \lim_{r\to 1^-} \int_{0}^{2\pi} \frac{1}{1-|\lambda|^2}\|f(\cdot, re^{i\theta})\|_{H^2(\mathbb D)}^2\frac{d\theta}{2\pi}\notag \\
    & =\frac{1}{1-|\lambda|^2}\|f\|_{H^2(\mathbb D^2)}^2.\label{h2-l-est}
\end{align}
Let $f\in\mathcal D(\mu_1,\mu_2).$ Now we derive an estimate for  $D_{\mu_2}(L_\lambda f)$ in terms of $D_{\mu_1,\mu_2}(f).$ For any  $z_2\in \mathbb D$, and  $s\in (|\lambda|,1)$  the function $g_{s,z_2}$ defined by
\begin{equation*}
g_{s,z_2}(w):=\partial_2f(sw,z_2),\quad w\in \mathbb D,
\end{equation*}
is holomorphic in a neighbourhood of $\overline{\mathbb D},$ and hence  $g_{s,z_2}\in H^2(\mathbb D).$ Using the reproducing property of $H^2(\mathbb D)$, for any $w\in \mathbb D,$ 
\begin{align}\label{g-s-w-evaluation}
|g_{s,z_2}(w)|^2\leqslant \frac{1}{1-|w|^2}\|g_{s,z_2}\|^2_{H^2(\mathbb D)}=\frac{1}{1-|w|^2} \int_0^{2\pi} |\partial_2f(se^{i\theta},z_2)|^2 \frac{d\theta}{2\pi}.
\end{align}
Substituting $w=\frac{\lambda}{s}$ in \eqref{g-s-w-evaluation} gives
\begin{equation}\label{par-point-eva}
|\partial_2f(\lambda,z_2)|^2\leqslant \frac{1}{1-|\lambda/s|^2} \int_0^{2\pi} |\partial_2f(se^{i\theta},z_2)|^2 \frac{d\theta}{2\pi}.
\end{equation}
Multiplying $P_{\mu_2}(z_2)$ on the both sides of \eqref{par-point-eva} and integrating over $\mathbb D$ gives 
\begin{align*}
 D_{\mu_2}(L_\lambda f)&=\int_\mathbb D |\partial_2(f(\lambda,z_2))|^2P_{\mu_2}(z_2)dA(z_2)\\
 &\leqslant  \frac{s^2}{s^2-|\lambda|^2}\int_\mathbb D  \int_0^{2\pi} |\partial_2f(se^{i\theta},z_2)|^2 P_{\mu_2}(z_2)\frac{d\theta}{2\pi}dA(z_2).
\end{align*}
Applying Fubini-Tonelli theorem gives
\begin{align*}
 D_{\mu_2}(L_\lambda f)\leqslant \frac{s^2}{s^2-|\lambda|^2}  \int_0^{2\pi}\int_\mathbb D |\partial_2f(se^{i\theta},z_2)|^2 P_{\mu_2}(z_2)dA(z_2)\frac{d\theta}{2\pi}.
\end{align*}
Taking limit $s\to 1^{-}$ yields 
\begin{equation}\label{d-mu-l-est}
D_{\mu_2}(L_\lambda f)\leqslant \frac{1}{1-|\lambda|^2} D^{(2)}_{\mu_1,\mu_2}(f).
\end{equation}
Combining \eqref{h2-l-est} and \eqref{d-mu-l-est} gives
\begin{align*}
\|L_\lambda f\|^2_{\mathcal D(\mu_2)}\leqslant \frac{1}{1-|\lambda|^2}\|f\|^2_{\mathcal D(\mu_1,\mu_2)}.
\end{align*}
Similar estimate can be derived for $R_\lambda f.$ Hence the result.
\end{proof}

Here are some immediate consequences of Lemma~\ref{slices-d-mu}.
\begin{corollary}\label{slices-in-d-mu-12}
Let $f\in \mathcal D(\mu_1,\mu_2)$ and $\lambda=(\lambda_1, \lambda_2)\in \mathbb D^2$.  Then the functions $f(\lambda_1,z_2)$ and $f(z_1,\lambda_2)$  belong to $\mathcal D(\mu_1,\mu_2).$
\end{corollary}
\begin{proof}
By \eqref{eqnhardy1}, $f(\lambda_1,\cdot), f(\cdot,\lambda_2) \in H^2(\mathbb D)$, and therefore  $f(\lambda_1,z_2)$ and $f(z_1,\lambda_2)$  belong to $H^2(\mathbb D^2).$ Lemma~\ref{slices-d-mu} gives $f(z_1,\lambda_2)\in \mathcal D(\mu_1)$ and  $f(\lambda_1,z_2)\in \mathcal D(\mu_2).$ Using the definition,
\begin{align*}
D_{\mu_1,\mu_2}(f(z_1,\lambda_2))=D_{\mu_1}(f(\cdot,\lambda_2))<\infty\quad \text{and}\quad 
D_{\mu_1,\mu_2}(f(\lambda_1,z_2))=D_{\mu_2}(f(\lambda_1,\cdot))<\infty,
\end{align*}
which completes the proof.
\end{proof}

One of the intrinsic properties of the local Dirichlet space $\mathcal D(\delta_{\zeta})$, $\zeta\in\mathbb T$, is that the boundary value 
\begin{equation}\label{eq:nontangential}
\text{$g^*(\zeta)=\lim_{r\to 1^-}g(r\zeta)$ exists for all $g\in \mathcal D(\delta_{\zeta})$,}
\end{equation}
 see  \cite[Proposition~2.1]{RS1991}. In the following, we show that this does not hold for $\mathcal D(\delta_{\eta},\delta_\zeta)$, where $(\eta,\zeta)\in\mathbb T^2$, that is, there exists a $f\in \mathcal D(\delta_{\eta},\delta_\zeta)$ for which $f^*(\eta, \zeta)$ does not exists. We begin with the following corollary which is an immediate consequence of Lemma~\ref{slices-d-mu}.

\begin{corollary}
Let $\zeta\in \mathbb T$, and $\mu\in M_+(\mathbb T)$. Suppose $f\in \mathcal D(\delta_{\zeta}, \mu)$. Then  $$f^*(\zeta, z_2):=\lim_{r\to 1^-}f(r\zeta,z_2)$$ exists for each $z_2\in\mathbb D$.
\end{corollary}
\begin{proof}
By Lemma~\ref{slices-d-mu}, for each $z_2\in \mathbb D$, the  function $f(\cdot, z_2)$ belongs to $\mathcal D(\delta_\zeta).$ Then by \eqref{eq:nontangential}, the result follows.
\end{proof}
Note the for any $f\in \mathcal O(\mathbb D^2)$ and finite nonnegative Borel measures $\mu_1,\mu_2$ on $\mathbb T$, we have 
\begin{align}\label{formula-D}
D_{\mu_1, \mu_2}(f) = \lim_{r \to 1^{-}} \int_{0}^{2\pi} D_{\mu_1}(f(\cdot, re^{i\theta})) \frac{d\theta}{2\pi} 
+ \lim_{r \to 1^{-}} \int_{0}^{2\pi} D_{\mu_2}(f(re^{i\theta}, \cdot))\frac{d\theta}{2\pi}.
\end{align}
\begin{proposition}
 Let $(\eta, \zeta)\in \mathbb T^2$.  There exists a function $f\in \mathcal D(\delta_\eta, \delta_\zeta)$ for which the boundary value  
 $$f^*(\eta, \zeta):=\lim_{r\to 1^-}f(r\eta,r\zeta)$$
  does not exist.
\end{proposition}
\begin{proof}
Consider the analytic Wiener algebra $$\mathcal W^+:=\{\sum_{k=0}^\infty a_kz^k\in\mathcal O(\mathbb D): \sum_{k=0}^\infty |a_k|<\infty\}.$$ Choose a function $g\in \mathcal D(\mathbb D)\backslash \mathcal W^+$ such that $g$ has nonnegative Taylor coefficients. Define $f\in \mathcal O(\mathbb D^2)$ by 
$$f(z_1, z_2)=g(\bar{\eta}\bar{\zeta} z_1z_2),~(z_1,z_2)\in\mathbb D^2.$$
Then 
\begin{align}\label{eq:d-eta-zeta}
    D_{\delta_\eta,\delta_\zeta}^{(1)}(f)=\lim_{r\to 1^-}\int_{0}^{2\pi}D_{\delta_\eta}(f(\cdot,re^{i\theta}))\frac{d\theta}{2\pi}=\lim_{r\to 1^-}\int_{0}^{2\pi}D_{\delta_\eta}(g(\bar{\eta}\bar{\zeta} re^{i\theta}w))\frac{d\theta}{2\pi}.
\end{align}
Note that, by the local Douglas formula \cite[Theorem 7.2.5]{EKMR2014},
\begin{align*}
    D_{\delta_\eta}(g(\bar{\eta}\bar{\zeta} re^{i\theta}w))&=\left\|\frac{g(\bar{\eta}\bar{\zeta} re^{i\theta}w)-g(\bar{\zeta} re^{i\theta})}{w-\eta}\right\|^2_{H^2(\mathbb D)}\\
    &=\left\|\frac{g(rw)-g(r\bar{\zeta}e^{i\theta})}{w-\bar{\zeta}e^{i\theta}}\right\|^2_{H^2(\mathbb D)} (\mbox{by rotation invariance of $\|\cdot\|_{H^2(\mathbb D)}$})\\
    &=D_{\delta_{\bar{\zeta}e^{i\theta}}}(g(rw)).
\end{align*}
Hence from \eqref{eq:d-eta-zeta}, we obtain
\begin{align*}
    D_{\delta_\eta,\delta_\zeta}^{(1)}(f)=\lim_{r\to 1^-}\int_{0}^{2\pi}D_{\delta_{\bar{\zeta}e^{i\theta}}}(g(rw))\frac{d\theta}{2\pi}=\lim_{r\to 1^-}D_{m}(g(rw))=D_{m}(g).
\end{align*}
Similarly, $D_{\delta_\eta,\delta_\zeta}^{(2)}(f)=D_{m}(g).$ Hence 
\begin{align*} 
    \|f\|^2_{\mathcal D(\delta_\eta, \delta_\zeta)}= \|f\|^2_{H^2(\mathbb D^2)}+2D_{m}(g)=\sum_{k=0}^\infty(1+2k)|\hat{g}(k)|^2<\infty.
\end{align*}
Since $g\in\mathcal D(\mathbb D)$, it follows that $f\in \mathcal D(\delta_\eta, \delta_\zeta).$ But $f(r\eta,r\zeta)=\sum_{k=0}^\infty \hat{g}(k) r^{2k}$, $0<r<1.$ If $f^*(\eta, \zeta)$ exists, then we must have $\sum_{k=0}^\infty \hat{g}(k)< \infty$, which is a contradiction 
as $g\notin \mathcal  W^+$ and $g$ has nonnegative Taylor coefficients. Hence $f^*(\eta, \zeta)$ does not exists.
\end{proof}
For $i=1,2,$ define the linear map $\Gamma_i:\mathcal O(\mathbb D) \rightarrow \mathcal O(\mathbb D^2)$ as 
\begin{equation*}
(\Gamma_i g)(z_1,z_2):=g(z_i), \quad (z_1,z_2)\in \mathbb D^2, g\in \mathcal O(\mathbb D).
\end{equation*}
Clearly, $\Gamma_i$ maps $\mathcal D(\mu_i)$ isometrically into  $\mathcal D(\mu_1,\mu_2)$ for $i=1,2.$

\begin{lemma}\label{one-to-two-mult}
For $i=1,2,$ $\Gamma_i$ maps $\mathrm{Mult}(\mathcal D(\mu_i))$ contractively into $\mathrm{Mult}(\mathcal D(\mu_1,\mu_2)).$
\end{lemma}
\begin{proof}
Without loss of generality, let $g\in \mathrm{Mult}(\mathcal D(\mu_1)).$ 
Thus we have $$g\in H^\infty(\mathbb D),~\mbox{and}~\|g\|_{H^\infty(\mathbb D)}\leqslant\|g\|_{\operatorname{Mult}(\mathcal D(\mu_1)}.$$ 
Let $f\in \mathcal D(\mu_1,\mu_2)$. Then   $f\in H^2(\mathbb D^2),$ and thus $(\Gamma_1 g)f\in H^2(\mathbb D^2)$. Moreover,  
\begin{align}
\|(\Gamma_1 g)f\|^2_{H^2(\mathbb D^2)}&\overset{\eqref{eqnhardy}}{=}\lim_{r\to1^-}\int_{0}^{2\pi} \|((\Gamma_1 g)f)(\cdot,re^{i\theta})\|^2_{H^2(\mathbb D)}\frac{d\theta}{2\pi} \notag\\& =\lim_{r\to1^-} \int_{0}^{2\pi} \|g(\cdot)f(\cdot,re^{i\theta})\|^2_{H^2(\mathbb D)}\frac{d\theta}{2\pi}.\label{eqnmult}
\end{align}
Thus 
\begin{eqnarray*}
\|(\Gamma_1 g)f\|^2_{\mathcal D(\mu_1,\mu_2)}&=&\left( \|(\Gamma_1 g)f\|_{H^2(\mathbb D^2)}^2+D_{\mu_1,\mu_2}^{(1)}((\Gamma_1 g)f)\right)+D_{\mu_1,\mu_2}^{(2)}((\Gamma_1 g)f)\\
& \overset{\eqref{formula-D}\& \eqref{eqnmult}}{\leqslant}&  \lim_{r\to1^-}\int_{0}^{2\pi}\left( \|g(\cdot)f(\cdot,re^{i\theta})\|^2_{H^2(\mathbb D)}+D_{\mu_1}(g(\cdot)f(\cdot,re^{i\theta}))\right)\frac{d\theta}{2\pi}\\
&&+\,\|g\|^2_{H^\infty(\mathbb D)}\lim_{r\to1^-}\int_{0}^{2\pi}\int_{\mathbb D}\left|\partial_2 f(re^{i\theta},z_2)\right|^2 P_{\mu_2}(z_2)\,d A(z_2)\,\frac{d\theta}{2\pi}\\
&=&\lim_{r\to1^-}\int_{0}^{2\pi} \|g(\cdot) f(\cdot,re^{i\theta})\|^2_{\mathcal D(\mu_1)}\frac{d\theta}{2\pi}+\|g\|^2_{H^\infty(\mathbb D)}D_{\mu_1,\mu_2}^{(2)}(f).
\end{eqnarray*}
Since $g\in \operatorname{Mult}(\mathcal D(\mu_1))$, by Lemma \ref{slices-d-mu},  
$$\|g(\cdot) f(\cdot,re^{i\theta})\|_{\mathcal D(\mu_1)}\leq \|g\|_{\mathrm{Mult}(\mathcal D(\mu_1))}  \|f(\cdot,re^{i\theta})\|_{\mathcal D(\mu_1)}.$$
Therefore
\begin{eqnarray*}\|(\Gamma_1 g)f\|^2_{\mathcal D(\mu_1,\mu_2)}
  &\leqslant& \|g\|^2_{\mathrm{Mult}(\mathcal D(\mu_1))} \lim_{r\to1^-}\int_{0}^{2\pi} \|f(\cdot,re^{i\theta})\|^2_{\mathcal D(\mu_1)}\frac{d\theta}{2\pi}+\|g\|^2_{H^\infty(\mathbb D)}D_{\mu_1,\mu_2}^{(2)}(f)\\
&=&\|g\|^2_{\mathrm{Mult}(\mathcal D(\mu_1))}\left(\|f\|^2_{H^2(\mathbb D^2)}+\lim_{r\to1^-}\int_{0}^{2\pi} D_{\mu_1}(f(\cdot,re^{i\theta}))\frac{d\theta}{2\pi}\right)\\
&&+\|g\|^2_{H^\infty(\mathbb D)}D_{\mu_1,\mu_2}^{(2)}(f)\\
&\leqslant&  \|g\|^2_{\mathrm{Mult}(\mathcal D(\mu_1))}\|f\|^2_{\mathcal D(\mu_1,\mu_2)}.
\end{eqnarray*}
 Hence $\Gamma_1$ is a contraction. Similarly, $\Gamma_2$ is also a contraction.
 \end{proof}
\begin{lemma}\label{slice-mult}
Let  $\varphi\in \mathrm{Mult}(\mathcal D(\mu_1,\mu_2))$ and $\lambda=(\lambda_1, \lambda_2)\in \mathbb D^2$.  Then the functions $\varphi(\lambda_1,z_2)$ and $\varphi(z_1,\lambda_2)$  belong to $\mathrm{Mult}(\mathcal D(\mu_1,\mu_2)).$
\end{lemma}
\begin{proof}
In view of Lemma~\ref{one-to-two-mult}, it suffices to show  $$L_{\lambda_1} \varphi \in  \mathrm{Mult}(\mathcal D(\mu_2))\text{ and }R_{\lambda_2} \varphi \in  \mathrm{Mult}(\mathcal D(\mu_1)).$$ 
Let $g\in \mathcal D(\mu_2)$. Since $\varphi\in \mathrm{Mult}(\mathcal D(\mu_1,\mu_2))$, it is easy to see that $\varphi\Gamma_2g \in \mathcal D(\mu_1,\mu_2).$ By Lemma~\ref{slices-d-mu},  we get  
$$\|L_{\lambda_1}(\varphi \Gamma_2g)\|^2_{\mathcal D(\mu_2)}\leqslant \frac{1}{1-|\lambda_1|^2}\|\varphi \Gamma_2g\|^2_{\mathcal D(\mu_1,\mu_2)}<\infty.$$
Note that $L_{\lambda_1}(\varphi \Gamma_2g)=(L_{\lambda_1}\varphi) g$. Indeed, for $w\in \mathbb D,$ 
$$L_{\lambda_1}(\varphi \Gamma_2g)(w)=(\varphi \Gamma_2g)(\lambda_1, w)=\varphi(\lambda_1, w)(\Gamma_2g)(\lambda_1, w)=(L_{\lambda_1}\varphi)(w)g(w)=((L_{\lambda_1}\varphi) g)(w).$$
Hence $(L_{\lambda_1}\varphi) g\in \mathcal D(\mu_2)$ for any $g\in \mathcal D(\mu_2).$ This proves that $L_{\lambda_1} \varphi \in  \mathrm{Mult}(\mathcal D(\mu_2)).$ Similarly, $R_{\lambda_2} \varphi\in \mathrm{Mult}(\mathcal D(\mu_1)).$ 
\end{proof}

\section{Gleason's problem and applications}\label{section-3}
In this section, we solve  Gleason's problem for both  the spaces $\mathcal D(\mu_1,\mu_2)$ and their  multiplier algebras. Using Corollary~\ref{slices-in-d-mu-12} and \cite[Lemma~5.3]{BCG2025}, we obtain a solution of Gleason's problem for  $\mathcal D(\mu_1,\mu_2)$. For the sake of completeness, we include the details here.
\begin{lemma}\label{gleason}
 Gleason's problem can be solved for $\mathcal D(\mu_1,\mu_2).$ 
\end{lemma}
\begin{proof} 
 Let $f\in \mathcal D(\mu_1,\mu_2)$ and $\lambda=(\lambda_1,\lambda_2)\in \mathbb D^2.$ By  Corollary~\ref{slices-in-d-mu-12},  $f(z_1,\lambda_2)$ and  $f(\lambda_1,z_2)$ belong to $\mathcal D(\mu_1,\mu_2).$  By Lemma \cite[Lemma 2.1]{BCJ2025}, there exists $f_1\in\mathcal O(\mathbb D^2)$ (may depends on $f$ and $\lambda$) such that 
\begin{equation}\label{gleason-1st}
f(z_1,z_2)-f(z_1,\lambda_2)=(z_2-\lambda_2)f_1(z_1,z_2),~(z_1,z_2)\in\mathbb D^2.
\end{equation}
Since $(z_2-\lambda_2)f_1(z_1,z_2)\in \mathcal D(\mu_1,\mu_2)$ and $\mathcal D(\mu_1,\mu_2)$ has the division property (see \cite[Theorem~2.2]{BCG2025}), it follows that $f_1\in \mathcal D(\mu_1,\mu_2).$
Similarly, there exists $f_2(z_1)\in \mathcal \mathcal D(\mu_1,\mu_2)$ (may depends on $f$ and $\lambda$) such that 
\begin{equation}\label{gleason-2nd}
f(z_1,\lambda_2)-f(\lambda_1,\lambda_2)=(z_1-\lambda_1)f_2(z_1),~z_1\in\mathbb D.
\end{equation}
Finally, adding \eqref{gleason-1st} and \eqref{gleason-2nd} yields the result.
\end{proof}

We recall some basic notions from multivariate spectral theory
(see \cite{Cu1981,T1970}). Let $T=(T_1,T_2)$ be a commuting pair of
bounded operators on a Hilbert space $\mathcal H$. The \emph{Koszul complex}
associated with $T$ is
$$K(T)\colon \{0\} \xrightarrow{0} \mathcal H \xrightarrow{B_2} \mathcal H \oplus \mathcal H
\xrightarrow{B_1} \mathcal H \xrightarrow{0} \{0\},$$
where the boundary maps are defined by
$$B_2(h)=(T_2h,-T_1h),\qquad B_1(h_1,h_2)=T_1h_1+T_2h_2.$$
Since $T_1$ and $T_2$ commute, $B_1B_2=T_1T_2-T_2T_1=0$.
For $k=0,1,2$, let $H^k(T)$ denote the $k$th cohomology group of $K(T)$.
The pair $T$ is \emph{Taylor-invertible} if $H^k(T)=\{0\}$ for every
$k=0,1,2$ \cite{T1970}. It is \emph{Fredholm} if each $H^k(T)$ is
finite-dimensional \cite{Cu1981}. For $\lambda=(\lambda_1,\lambda_2)
\in\mathbb C^2$, write $T-\lambda=(T_1-\lambda_1 I,T_2-\lambda_2 I)$.
The \emph{Taylor spectrum} and the \emph{essential spectrum} of $T$
are then given by
\begin{align*}
\sigma(T)&= \{\lambda\in\mathbb C^2 : T-\lambda \text{ is not Taylor-invertible}\}, \\
\sigma_e(T)&= \{\lambda\in\mathbb C^2 : T-\lambda \text{ is not Fredholm}\}.
\end{align*}
For a Fredholm pair $T$, the \emph{Fredholm index} is the Euler
characteristic of its Koszul complex,
$$\operatorname{ind}(T) = \dim H^0(T)-\dim H^1(T)+\dim H^2(T).$$

\begin{theorem}\label{essential-spectrum}
The Taylor essential spectrum $\sigma_e(M_z)$ of $M_z$ on $\mathcal D(\mu_1,\mu_2)$ is the topological boundary $\partial \mathbb D^2$ of $\mathbb D^2.$ 
\end{theorem}
\begin{proof}
We already know that $\sigma(M_z)=\overline{\mathbb D}^2$ (see \cite[Eq.~(7.1)]{BCG2025}). Since $\sigma_e(M_z)\subseteq \sigma(M_z),$  so $\sigma_e(M_z)\subseteq \overline{\mathbb D}^2.$ Consider the Koszul complex $K(M_z-\lambda)$. 
Then for $f,f_1,f_2\in \mathcal D(\mu_1,\mu_2),$
\begin{align*} B_2(f):=((z_2-\lambda_2)f,-(z_1-\lambda_1)f), \, B_1(f_1,f_2):= (z_1-\lambda_1)f_1 +(z_2-\lambda_2)f_2.
\end{align*}
Now let  $\lambda\in \mathbb D^2.$ Note that, $\ker B_2=\{0\},$ and by \cite[Corollary~4.6]{BCG2025}, $\mathrm{ran}\,B_2=\ker B_1.$ Also, by Lemma~\ref{gleason}, it follows that $\dim(\mathcal D(\mu_1,\mu_2)/\mathrm{ran}\,B_1)=1.$  Thus for each $\lambda\in \mathbb  D^2,$ $M_z-\lambda$ is Fredholm with $\mathrm{ind}\,(M_z-\lambda)=1$ and hence the Taylor essential spectrum $\sigma_e(M_z)\subseteq \partial \mathbb D^2.$ 
Moreover, for each $\lambda\in \mathbb C^2\backslash \overline{\mathbb D}^2,$ $(M_z-\lambda)$ is invertible so $\mathrm{ind}\,(M_z-\lambda)=0.$ If possible, let $\lambda_0\in \partial \mathbb D^2\backslash \sigma_e(M_z),$ then we can find two different sequences one from $\mathbb D^2$ and another from $\mathbb C^2\backslash \overline{\mathbb D}^2$ converging to $\lambda_0.$ Since the Fredholm index is continuous, we get a contradiction. Hence $\sigma_e(M_z)=\partial \mathbb D^2.$ 
\end{proof}
Since Gleason's problem has a solution for each $\lambda\in \mathbb D^2,$ so  \cite[Corollary~3.9 \& Lemma~5.1]{BCG2025} imply $M_z^*$ belongs to the Cowen-Douglas class $\boldsymbol B_1(\mathbb D^2)$ of rank $1$ over $\mathbb D^2.$ This improves the statement $\mathrm{(ii)}$ of \cite[Theorem~2.4]{BCG2025} as follows:
\begin{theorem}
Let $T=(T_1, T_2)$ be a commuting pair on $\mathcal H.$ Then the following statements are equivalent:
\begin{enumerate}
\item[$\mathrm{(i)}$] $T$ is a cyclic analytic toral $2$-isometry with cyclic vector $f_0 \in \ker T^*=\ker T_1^*\cap \ker T_2^*$ and $\ker T^*$ is a wandering subspace for $T,$
\item[$\mathrm{(ii)}$] $T$ is a cyclic toral $2$-isometry with cyclic vector $f_0 \in \ker T^*$, $T^*$ belongs to ${\bf B}_1(\mathbb D^2)$ and $\ker T^*$ is a wandering subspace for~$T,$
\item[$\mathrm{(iii)}$] there exist $\mu_1, \mu_2 \in M_+(\mathbb T)$ such that $T$ is unitarily equivalent to $M_z$ on $\mathcal D(\mu_1, \mu_2).$
\end{enumerate}
\end{theorem}
\begin{proof}[Proof of Theorem~\ref{thm:gleason}] The solvability of Gleason's problem for the space $\mathcal D(\mu_1,\mu_2)$ follows from  Lemma~\ref{gleason}. To solve Gleason's problem for $\mathrm{Mult}(\mathcal D(\mu_1,\mu_2)),$ let $\varphi$ be a multiplier of $\mathcal D(\mu_1,\mu_2)$ and $(\lambda_1,\lambda_2)\in \mathbb D^2.$ By Lemma~\ref{slice-mult}, $\varphi(\lambda_1,z_2)$ belongs to $\mathrm{Mult}(\mathcal D(\mu_1,\mu_2)).$  By Lemma \cite[Lemma 2.1]{BCJ2025}, there exists $\varphi_1\in\mathcal O(\mathbb D^2)$ (may depends on $\varphi$ and $\lambda$) such that 
\begin{equation}\label{division-1}
\varphi(z_1,z_2)-\varphi(\lambda_1,z_2)=(z_1-\lambda_1)\varphi_1(z_1,z_2).
\end{equation}
 Now we show that $\varphi_1(z_1,z_2)$ is a multiplier of $\mathcal D(\mu_1,\mu_2).$ To see this, let $f\in \mathcal D(\mu_1,\mu_2),$ then 
\begin{align*}
(\varphi_1 f)(z_1,z_2)=\frac{(z_1-\lambda_1)\varphi_1(z_1,z_2)}{z_1-\lambda_1}f(z_1,z_2)\overset{\eqref{division-1}}{=}\frac{\varphi(z_1,z_2)f(z_1,z_2)-\varphi(\lambda_1,z_2)f(z_1,z_2)}{z_1-\lambda_1}.
\end{align*}
Note that the numerator belongs to $\mathcal D(\mu_1,\mu_2)$ and vanishes at $(\lambda_1,z_2).$ Since $\mathcal D(\mu_1,\mu_2)$ has division property (see \cite[Theorem~2.2]{BCG2025}), there exists $g\in \mathcal D(\mu_1,\mu_2)$ such that $$\varphi(z_1,z_2)f(z_1,z_2)-\varphi(\lambda_1,z_2)f(z_1,z_2)=(z_1-\lambda_1)g(z_1,z_2).$$
This shows $\varphi_1 f=g\in \mathcal D(\mu_1,\mu_2)$ for any $f\in  \mathcal D(\mu_1,\mu_2).$ Hence $\varphi_1$ is a multiplier of $\mathcal D(\mu_1,\mu_2).$ Similarly, there exists $\varphi_2\in \mathcal O(\mathbb D)$ (may depends on $\varphi$ and $\lambda$) such that 
\begin{equation}\label{division-2}
\varphi(\lambda_1,z_2)-\varphi(\lambda_1,\lambda_2)=(z_2-\lambda_2)\varphi_2(z_2).
\end{equation}
Arguing as above one can show that $\varphi_2(z_2)$ is a multiplier of $\mathcal D(\mu_1,\mu_2).$ Finally, adding \eqref{division-1} and \eqref{division-2} yields the result. 
\end{proof}

 We are now ready to prove Theorem~\ref{codim-thm}
\begin{proof}[Proof of Theorem~\ref{codim-thm}]
We begin by showing that 
\begin{align}\label{codim-one}
\dim \mathcal S/\left( (z_1-\lambda_1)\mathcal S+(z_2-\lambda_2)\mathcal S\right) =1
\end{align}
for all $\lambda\in \mathbb D^2\backslash Z(\varphi).$ To see this, let us fix $\lambda\in \mathbb D^2\backslash Z(\varphi).$ Without loss of generality, assume that $\varphi(\lambda_1,\lambda_2)=1.$ Since Gleason's problem has a solution for $\mathrm{Mult}(\mathcal D(\mu_1,\mu_2))$ (see Theorem \ref{thm:gleason}), there exist $\varphi_1,\varphi_2\in \mathrm{Mult}(\mathcal D(\mu_1,\mu_2))$ such that 
\begin{equation*}
\varphi(z_1,z_2)-\varphi(\lambda_1,\lambda_2)=(z_1-\lambda_1)\varphi_1(z_1,z_2)+(z_2-\lambda_2)\varphi_2(z_1,z_2), ~(z_1,z_2)\in\mathbb D^2.
\end{equation*}
This gives
\begin{equation*}
1=\varphi(z_1,z_2)-(z_1-\lambda_1)\varphi_1(z_1,z_2)-(z_2-\lambda_2)\varphi_2(z_1,z_2).
\end{equation*}
Let $f\in \mathcal S.$ Then
\begin{eqnarray}\label{f-phi-gleason}
f=f\cdot 1= f\left[ \varphi-(z_1-\lambda_1)\varphi_1-(z_2-\lambda_2)\varphi_2\right]=\varphi f-(z_1-\lambda_1)\varphi_1 f-(z_2-\lambda_2)\varphi_2 f.
\end{eqnarray}
Since Gleason's problem also has a solution for $\mathcal D(\mu_1,\mu_2)$ (see Lemma~\ref{gleason}), for $f$ there exist $f_1,f_2\in \mathcal D(\mu_1,\mu_2)$ (depend on $f$ and $\lambda$) such that 
\begin{equation}\label{f-gleason}
f(z_1,z_2)=f(\lambda_1,\lambda_2)+(z_1-\lambda_1)f_1+(z_2-\lambda_2)f_2.
\end{equation}
Combining \eqref{f-phi-gleason} and \eqref{f-gleason}  yields 
\begin{equation}\label{phi-codim-1}
f= f(\lambda_1,\lambda_2)\varphi+(z_1-\lambda_1)(\varphi f_1-\varphi_1 f)+(z_2-\lambda_2)(\varphi f_2-\varphi_2 f).
\end{equation}
We now claim that 
\begin{equation}\label{eqnclaim}
    \varphi \mathcal D(\mu_1,\mu_2)\subset S.
\end{equation}
To see this, let $g\in \mathcal D(\mu_1,\mu_2)$.
Let $\{p_n^{(i)}\}$ be a sequence of polynomials converging to $g$ in $\mathcal D(\mu_1,\mu_2)$. Since $\varphi$ is a multiplier of $\mathcal D(\mu_1,\mu_2)$, $\{\varphi p_n^{(i)}\}$ converges to $\varphi g.$ Since $\mathcal S$ is a submodule and $\varphi\in\mathcal S$, we get $\varphi p_n^{(i)}\in \mathcal S$ and hence $\varphi g\in \mathcal S$, establishing the claim. Consequently, we have $\varphi f_i\in \mathcal S$ for $i=1,2$. Moreover, since $\mathcal S$ is multiplier invariant, $\varphi_i f\in \mathcal S$ for $i=1,2.$ Therefore, by \eqref{phi-codim-1}, we conclude that $$f-f(\lambda_1,\lambda_2)\varphi\in (z_1-\lambda_1)\mathcal S+(z_2-\lambda_2)\mathcal S $$
establishing \eqref{codim-one}.

For the remaining part we apply \cite[Theorem~3.4]{GRS2005} (with $\mathscr E=\mathscr D=\mathbb C$). Note that $1\in \mathcal D(\mu_1,\mu_2)$ and $z_1,z_2$ are multipliers of $\mathcal D(\mu_1,\mu_2).$ Theorem~\ref{thm:gleason} shows that one can solve Gleason's problem for $\mathrm{Mult}(\mathcal D(\mu_1,\mu_2)).$ For $\lambda=(\lambda_1,\lambda_2)\in \mathbb D^2,$ consider the augmented Koszul complex $K(M_z-\lambda,\mathbb C)$ 
$$K(M_z-\lambda,\mathbb C):\{0\}\overset{0}{\longrightarrow} \mathcal H\overset{B_2}{\longrightarrow} \mathcal H \oplus \mathcal H\overset{B_1}{\longrightarrow}\mathcal H\overset{E_\lambda}{\longrightarrow} \mathbb C\overset{0}{\longrightarrow} \{0\},$$
where $\mathcal H=\mathcal D(\mu_1,\mu_2),$ 
\begin{align*} B_2(f):=((z_2-\lambda_2)f,-(z_1-\lambda_1)f), \; B_1(f_1,f_2):= (z_1-\lambda_1)f_1 +(z_2-\lambda_2)f_2, 
\end{align*}
\begin{equation*}
E_\lambda (f):=f(\lambda_1,\lambda_2),\quad f,f_1,f_2\in \mathcal D(\mu_1,\mu_2).
\end{equation*}
Note that the exactness at all stages except the last one is already established in Theorem~ \ref{essential-spectrum}. Moreover, the exactness at the last stage follows from  Lemma~\ref{gleason}. Also from \eqref{eqnclaim} we have $\varphi\mathcal D(\mu_1,\mu_2) \subset \mathcal S$. This completes the verification of all the hypotheses of \cite[Theorem~3.4]{GRS2005}. Hence the result follows. 
\end{proof}

\section{Left-inverse commuting submodules of \texorpdfstring{$\mathcal D(\mu_1,\mu_2)$}{}}\label{section-4}
Let $f\in H^2(\mathbb D^2)$ with the Taylor series representation $$f(z_1,z_2)=\sum_{k=0}^\infty\sum_{\ell=0}^\infty \hat{f}(k,\ell)z_1^k z_2^\ell.$$
Following \cite[p.~382]{Y1999}, for each $k,\ell\geqslant 0,$ define 
\begin{equation}\label{hardy-1-2-slices}
    f^{(1)}_\ell(w):=\sum_{n=0}^\infty\hat{f}(n,\ell)w^n, \quad f^{(2)}_k(w):=\sum_{n=0}^\infty\hat{f}(k,n)w^n,\quad w\in \mathbb D.
\end{equation}
Then each $f^{(1)}_\ell, f^{(2)}_k$ belongs to $H^2(\mathbb D).$ Moreover,
\begin{equation*}
    f(z_1,z_2)=\sum_{k=0}^\infty f^{(2)}_k(z_2) z_1^k=\sum_{\ell=0}^\infty f^{(1)}_\ell(z_1) z_2^\ell,
\end{equation*}
and
\begin{equation*}
\|f\|_{H^2(\mathbb D^2)}^2=\sum_{k=0}^\infty \|f^{(2)}_k\|^2_{H^2(\mathbb D)}=\sum_{\ell=0}^\infty \|f^{(1)}_\ell\|^2_{H^2(\mathbb D)}.
\end{equation*}
Let us recall the following lemma from \cite{RS1991}.
\begin{lemma}\rm(\cite[Corollary~2.3]{RS1991}\rm)\label{lemboundary}
Let $g\in \mathcal D(\mu)$. Then the nontangential boundary value $g^*(\zeta)=\lim_{r\to 1^-}g(r\zeta)$ exists $\mu$-a.e. on $\mathbb T$ and 
\begin{equation}\label{one-variable-richter-formula}
\|wg\|_{\mathcal D(\mu)}^2-\|g\|_{\mathcal D(\mu)}^2=\int_{\mathbb T}|g^*(\zeta)|^2 d\mu(\zeta),
\end{equation}
where  $w$ denotes the coordinate function on $\mathcal D(\mu).$    
\end{lemma}
\begin{lemma}\label{weight-functions}
    Let $f\in \mathcal D(\mu_1,\mu_2).$ Then  the following holds:
    \begin{enumerate}
        \item[\rm(i)]  for each $k,\ell \geqslant 0$, $f^{(1)}_\ell\in \mathcal D(\mu_1)$, $f^{(2)}_k\in \mathcal D(\mu_2),$ and
        \begin{equation*}
     \sum_{\ell=0}^\infty D_{\mu_1}(f^{(1)}_\ell)=D^{(1)}_{\mu_1, \mu_2}(f),~\sum_{k=0}^\infty D_{\mu_2}(f^{(2)}_k)=D^{(2)}_{\mu_1, \mu_2}(f),
\end{equation*}
        \item[\rm(ii)] we have $\sum_{\ell=0}^\infty |f_\ell^{(1)*}|^2\in L^1(\mu_1)$,   $\sum_{k=0}^\infty |f_k^{(2)*}|^2\in L^1(\mu_2),$ and
         \begin{eqnarray}
    \|z_1f\|^2_{\mathcal D(\mu_1, \mu_2)}-\|f\|^2_{\mathcal D(\mu_1, \mu_2)}&=& \int_\mathbb T \sum_{\ell=0}^\infty |f_\ell^{(1)*}(\zeta)|^2 d\mu_1(\zeta), \label{z-1-f-minus-f}\\
   \|z_2f\|^2_{\mathcal D(\mu_1, \mu_2)}-\|f\|^2_{\mathcal D(\mu_1, \mu_2)}&=& \int_\mathbb T \sum_{k=0}^\infty |f_k^{(2)*}(\zeta)|^2 d\mu_2(\zeta)  \label{z-2-f-minus-f}, 
\end{eqnarray}
where $f_\ell^{(1)*}, f_k^{(2)*}$ respectively denote the nontangential limits of $f^{(1)}_\ell, f^{(2)}_k.$ 
    \end{enumerate} 
\end{lemma}
\begin{proof}
    Let $f\in \mathcal D(\mu_1,\mu_2).$ Applying Fubini-Tonelli theorem, we obtain  
\begin{eqnarray*}
    D^{(1)}_{\mu_1, \mu_2}(f)=\lim_{r\to 1^-} \int_{\mathbb D} \left( \int_{0}^{2\pi} |\partial_1 f(z_1, re^{i \theta})|^2 \frac{d\theta}{2\pi}\right) P_{\mu_1}(z_1) \,dA(z_1).
\end{eqnarray*}
For each $z_1\in \mathbb D$ and $r\in (0,1)$ define 
\begin{equation*}
    g_{r,z_1}(e^{i\theta}):=\partial_1f(z_1,re^{i\theta})=\sum_{\ell=0}^\infty f^{(1)'}_\ell(z_1)r^\ell e^{i\ell\theta}.
\end{equation*}
  Note that for each $z_1\in \mathbb D,$  $\partial_1f(z_1,\cdot)$ is holomorphic on $\mathbb D$ so for each $r\in(0,1),$ $g_{r,z_1}\in L^\infty(\mathbb T).$ Hence by Parseval's identity
    \begin{align*}
      \int_{0}^{2\pi} |\partial_1 f(z_1, re^{i \theta})|^2 \frac{d\theta}{2\pi}=\int_{0}^{2\pi}|g_{r,z_1}(e^{i\theta})|^2\frac{d\theta}{2\pi}= \sum_{\ell=0}^\infty |f^{(1)'}_\ell(z_1)|^2r^{2\ell}.
    \end{align*}
   An application of monotone convergence theorem (see \cite[Corollary~1.10]{SS2005}) gives
     \begin{align*}
    \int_{\mathbb D}\int_{0}^{2\pi} |\partial_1 f(z_1, re^{i \theta})|^2 \frac{d\theta}{2\pi}P_{\mu_1}(z_1) \,dA(z_1)
    &=\int_{\mathbb D}\sum_{\ell=0}^\infty r^{2\ell}|f^{(1)'}_\ell(z_1)|^2 P_{\mu_1}(z_1) \,dA(z_1)\\
    &= \sum_{\ell=0}^\infty r^{2\ell}\int_{\mathbb D}|f^{(1)'}_\ell(z_1)|^2 P_{\mu_1}(z_1) \,dA(z_1)\\
    &=\sum_{\ell=0}^\infty r^{2\ell} D_{\mu_1}(f^{(1)}_\ell).
\end{align*} 
Taking $\lim_{r\to 1^-}$ yields 
\begin{equation}\label{d-mu-1-f-n}
     \sum_{\ell=0}^\infty D_{\mu_1}(f^{(1)}_\ell)=D^{(1)}_{\mu_1, \mu_2}(f)<\infty.
\end{equation}
This in particular, implies $f^{(1)}_\ell\in \mathcal D(\mu_1)$ for each $\ell\geqslant 0.$ Replacing $f$ by $z_1f$ in \eqref{d-mu-1-z-1-f-n} we get
\begin{equation}\label{d-mu-1-z-1-f-n}
     \sum_{\ell=0}^\infty D_{\mu_1}(wf^{(1)}_\ell)=D^{(1)}_{\mu_1, \mu_2}(z_1f) <\infty,
\end{equation}
Subtracting \eqref{d-mu-1-f-n} from \eqref{d-mu-1-z-1-f-n} and applying Lemma~\ref{lemboundary} gives 
\begin{align*}
    D^{(1)}_{\mu_1, \mu_2}(z_1f)-D^{(1)}_{\mu_1, \mu_2}(f)= \sum_{\ell=0}^\infty D_{\mu_1}(wf^{(1)}_\ell)-\sum_{\ell=0}^\infty D_{\mu_1}(f^{(1)}_\ell)&= \sum_{\ell=0}^\infty \left (D_{\mu_1}(wf^{(1)}_\ell)-D_{\mu_1}(f^{(1)}_\ell)\right)\\
    &= \sum_{\ell=0}^\infty \int_\mathbb T |f^{(1)*}_\ell(\zeta)|^2 d\mu_1(\zeta).
\end{align*}
Again applying \cite[Corollary~1.10]{SS2005} we get 
\begin{align}\label{d-1-z-1-f-minus-f}
    D^{(1)}_{\mu_1, \mu_2}(z_1f)-D^{(1)}_{\mu_1, \mu_2}(f)= \int_\mathbb T \sum_{\ell=0}^\infty |f^{(1)*}_\ell(\zeta)|^2 d\mu_1(\zeta),
\end{align}
 and hence $\sum_{\ell=0}^\infty |f^{(1)*}_\ell|^2\in L^1(\mu_1).$ 
 Combining \eqref{d-1-z-1-f-minus-f} with the facts that $\|z_1f\|_{H^2(\mathbb D^2)}=\|f\|_{H^2(\mathbb D^2)}$ and $D^{(2)}_{\mu_1, \mu_2}(z_1f)=D^{(2)}_{\mu_1, \mu_2}(f)$ establishes \eqref{z-1-f-minus-f}. By arguing similarly, one can show 
 that $f^{(2)}_k\in \mathcal D(\mu_2)$ for each $k\geqslant 0,$  $\sum_{k=0}^\infty |f^{(2)*}_k|^2\in L^1(\mu_2)$ and \eqref{z-2-f-minus-f} holds.
 \end{proof}
\begin{proposition}\label{extremal-cyclic-subspace}
    Let $\varphi$ be an extremal function in $\mathcal D(\mu_1,\mu_2).$ Then the  cyclic subspace 
    \begin{align*}
[\varphi]=\overline{\mathrm{span}}\{z^\alpha \varphi: \alpha\in \mathbb Z_+^2\}
    \end{align*}
    is of the form $\varphi \mathcal D(\nu_1,\nu_2),$ where 
    \begin{eqnarray*}
      d\nu_1(\zeta)=\sum_{\ell=0}^\infty |\varphi^{(1)*}_\ell(\zeta)|^2 d\mu_1(\zeta),\quad 
        d\nu_2(\zeta)=\sum_{k=0}^\infty |\varphi^{(2)*}_k(\zeta)|^2 d\mu_2(\zeta).
    \end{eqnarray*}
     Moreover, $M_\varphi$ is an isometry from $\mathcal D(\nu_1,\nu_2)$ to $\mathcal D(\mu_1,\mu_2).$
\end{proposition}
\begin{proof}
Let $$\mathcal S=[\varphi]=\overline{\mathrm{span}}\{z^\alpha \varphi: \alpha\in \mathbb Z_+^2\}$$
and define $$T=(T_1, T_2),~~\text{where}~~T_i=M_{z_i}|_\mathcal S,\qquad i=1,2.$$
Since $\varphi$ is extremal, we have 
\begin{align*}
      \langle z_1^k\varphi, z_1^pz_2^{q+1}\varphi\rangle_{\mathcal D(\mu_1,\mu_2)}=0,~
            \langle z_2^\ell\varphi, z_1^{p+1}z_2^q\varphi\rangle_{\mathcal D(\mu_1,\mu_2)}=0,\quad  k,\ell,p,q\geqslant 0.
\end{align*}
Hence for any $\alpha\in \mathbb Z_+^2$, 
$$\left\langle T_i^*\varphi,z^\alpha\varphi\right\rangle_{\mathcal S}=\left\langle \varphi,
z^{\alpha+\epsilon_i}\varphi\right\rangle_{\mathcal S}=0,~i=1,2.$$
Thus $\varphi \in \ker T^*=\ker T_1^*\cap\ker T_2^*.$ Since $T$ is cyclic, $\ker T^*=\mathbb C \varphi.$ 
 Note that $(T_1,T_2)$ is a cyclic toral $2$-isometry with cyclic vector $\varphi\in \ker T^*$ and $\ker T^*=\mathbb C\varphi$ is a wandering subspace for $T$ (see \cite[Definition~1.3]{BCG2025}). So by \cite[Theorem~2.4]{BCG2025}, there exist unique finite nonnegative Borel measures $\nu_1,\nu_2$ on $\mathbb T,$  and a unitary operator 
\begin{align*}
U:\mathcal D(\nu_1,\nu_2)\to \mathcal S
\end{align*}
such that $UM_{z_i}=T_iU$ for $i=1,2.$ Multiplying $U$ with a unimodular constant, we may assume that $U1=\varphi.$ Using the intertwining relation, for any polynomial $p$ we get $Up=p\varphi.$ 
By standard arguments using polynomial density and reproducing property, it follows that $U=M_{\varphi},$
$\mathcal S=\varphi \mathcal D(\nu_1,\nu_2)$ and 
\begin{equation*}
M_\varphi:\mathcal D(\nu_1,\nu_2)\rightarrow \mathcal D(\mu_1,\mu_2)\; \text{ is an isometry}.
\end{equation*}
Next we derive explicit forms of the measures $\nu_1$ and $\nu_2.$ 
We claim that $\nu_1$ and $\nu_2$ are of the form 
$$ d\nu_1(\zeta)=\left(\sum_{\ell=0}^{\infty}\left|\varphi_\ell^{(1)*}(\zeta)\right|^2\right)d\mu_1(\zeta),\; \,d\nu_2(\zeta)=\left(\sum_{k=0}^{\infty}\left|\varphi_k^{(2)*}(\zeta)\right|^2\right)d\mu_2(\zeta),$$
where $\varphi_\ell^{(1)}, \varphi_k^{(2)}$ are defined by \eqref{hardy-1-2-slices}.
By Lemma~\ref{weight-functions}, $\nu_1$ and $\nu_2$ define finite nonnegative  Borel measures on $\mathbb T.$ 

Let $p(w)$ be a one-variable polynomial. Then
\begin{eqnarray*}
    \|p\|^2_{\mathcal D(\nu_1)}&=&\|p(z_1)\|^2_{\mathcal D(\nu_1,\nu_2)}{=}\|p(z_1)\varphi\|^2_{\mathcal D(\mu_1,\mu_2)},\\
     \|wp\|^2_{\mathcal D(\nu_1)}&=&\|z_1p(z_1)\|^2_{\mathcal D(\nu_1,\nu_2)}=\|z_1p(z_1)\varphi\|^2_{\mathcal D(\mu_1,\mu_2)}. 
\end{eqnarray*}
Hence
\begin{align*}
\|wp\|^2_{\mathcal D(\nu_1)}-\|p\|^2_{\mathcal D(\nu_1)}=\|z_1p(z_1)\varphi\|^2_{\mathcal D(\mu_1,\mu_2)}-\|p(z_1)\varphi\|^2_{\mathcal D(\mu_1,\mu_2)}.
\end{align*}
Since  $$p(z_1)\varphi(z_1,z_2)=\sum_{\ell=0}^{\infty}\left(p(z_1)\varphi_\ell^{(1)}(z_1)\right)z_2^\ell,$$
it follows from \eqref{one-variable-richter-formula} and \eqref{z-1-f-minus-f} that 
\begin{eqnarray*}
    \int_{\mathbb T}|p(\zeta)|^2\,d\nu_1(\zeta)=\int_{\mathbb T}\sum_{\ell=0}^\infty |(p\varphi)_\ell^{(1)*}(\zeta)|^2\,d\mu_1(\zeta)=\int_{\mathbb T}|p(\zeta)|^2\left (\sum_{\ell=0}^\infty |\varphi_\ell^{(1)*}(\zeta)|^2\right)\,d\mu_1(\zeta).
\end{eqnarray*}
Since this is true for any one-variable polynomial $p$, combining the polarization identity with the uniqueness of the trigonometric moment problem, we get $$d\nu_1(\zeta)=\sum_{\ell=0}^\infty |\varphi^{(1)*}_\ell(\zeta)|^2 d\mu_1(\zeta).$$ 
The formula for $\nu_2$ can be derived similarly.
\end{proof}

Let $\kappa_\lambda^\mu$ denotes the reproducing kernel of $\mathcal D(\mu).$ Then from \cite[Corollaries~2~\&~3]{Sh2001-2} we know that $\kappa^{\mu_1}_\lambda$ is a multiplier of $\mathcal D(\mu_1).$ In fact,
   \begin{equation*}
    |\kappa^{\mu_1}(w,\lambda)|\leqslant \frac{2}{1-|\lambda|},\quad \text{for all } w\in \mathbb D.
    \end{equation*}
Moreover, \cite[Theorem~2]{Sh2001-2} shows that 
 \begin{equation}\label{2re-k-1}
 2\mathrm{Re}\kappa^{\mu_1}(w,\lambda)-1\geqslant 0,\quad \text{for all } w\in \mathbb D.
 \end{equation}
Also from \cite[Corollary~4 \& Eq.~(14)]{Sh2001-2}, for any $g,h\in \mathcal D(\mu),$ we have
\begin{align}\label{eq:shimorin-pq}
    \langle g \kappa_\lambda^\mu,h\rangle_{\mathcal D(\mu)}+  \langle g ,h\kappa_\lambda^\mu\rangle_{\mathcal D(\mu)}-\langle g ,h\rangle_{\mathcal D(\mu)}&=\int_\mathbb T g^*(\zeta)\overline{h^*(\zeta)}P_\lambda(\zeta) dm(\zeta) \notag\\
    &~~+ \int_\mathbb T (2\mbox{Re}\kappa^\mu_\lambda(\zeta)-1)D_{\delta_\zeta}(g,h)d\mu(\zeta),
\end{align}
where $D_{\delta_\zeta}(g,h)$ denotes the Dirichlet semi inner product of $g$ and $h$.

\begin{lemma}\label{lemma:k-mu-mult}
    Let $f\in \mathcal D(\mu_1,\mu_2)$ and $\lambda\in \mathbb D.$ Then 
    \begin{align*}
\|L_{\lambda}(f)\|^2_{H^2(\mathbb D)}=\sum_{\ell=0}^\infty |f^{(1)}_\ell(\lambda)|^2 &\leqslant  2\sum_{\ell=0}^\infty  \mathrm{Re}\langle f^{(1)}_\ell \kappa^{\mu_1}_{\lambda}, f^{(1)}_\ell\rangle_{\mathcal D(\mu_1)}- \sum_{\ell=0}^\infty \|f^{(1)}_\ell\|^2_{\mathcal D(\mu_1)},\\ 
\|R_{\lambda}(f)\|^2_{H^2(\mathbb D)}=\sum_{k=0}^\infty |f^{(2)}_k(\lambda)|^2 &\leqslant  2\sum_{k=0}^\infty  \mathrm{Re}\langle f^{(2)}_k \kappa^{\mu_2}_{\lambda}, f^{(2)}_k\rangle_{\mathcal D(\mu_2)}- \sum_{k=0}^\infty \|f^{(2)}_k\|^2_{\mathcal D(\mu_1)},
\end{align*}
where $\kappa_{\lambda}^{\mu_i}$ denotes the reproducing kernel of $\mathcal D(\mu_i),$ $i=1,2.$ 
\end{lemma}
\begin{proof}
Let $f\in \mathcal D(\mu_1,\mu_2).$ Then by Lemma \ref{weight-functions}, for each $n\geqslant 0,$ $f^{(1)}_\ell\in \mathcal D(\mu_1),$ and   $$D_{\mu_1,\mu_2}^{(1)}(f)=\sum_{\ell=0}D_{\mu_1}(f_\ell^{(1)}).$$
Since $\kappa^{\mu_1}_\lambda\in \mathrm{Mult}(\mathcal D(\mu_1)),$ by Lemma~\ref{one-to-two-mult}, $\kappa^{\mu_1}_\lambda(z_1)$ belongs to  $\mathrm{Mult}(\mathcal D(\mu_1,\mu_2)).$ 
Substituting $g,h=f_\ell^{(1)}$ in \eqref{eq:shimorin-pq} gives
\begin{align*}
&\langle f^{(1)}_\ell \kappa^{\mu_1}_\lambda, f^{(1)}_\ell\rangle_{\mathcal D(\mu_1)}+\langle f^{(1)}_\ell , f^{(1)}_\ell \kappa^{\mu_1}_\lambda\rangle_{\mathcal D(\mu_1)}-\langle f^{(1)}_\ell , f^{(1)}_\ell \rangle_{\mathcal D(\mu_1)}\\
&= \int_{\mathbb T} |f^{(1)*}_\ell(\zeta)|^2 P_\lambda(\zeta) dm(\zeta)+ \int_{\mathbb T} (2\mathrm{Re}\kappa^{\mu_1}_\lambda(\zeta)-1)D_{\delta_\zeta}(f^{(1)}_\ell)d\mu_1(\zeta).
\end{align*}
Since $2\mathrm{Re}\kappa^{\mu_1}(w,\lambda)-1\geqslant 0$ for all $w\in \mathbb D$ (see \eqref{2re-k-1}), we get
\begin{align*}
2 \mathrm{Re}\langle f^{(1)}_\ell \kappa^{\mu_1}_\lambda, f^{(1)}_\ell\rangle_{\mathcal D(\mu_1)}-\|f^{(1)}_\ell\|^2_{\mathcal D(\mu_1)}\geqslant \int_{\mathbb T} |f^{(1)*}_\ell(\zeta)|^2 P_\lambda(\zeta) dm(\zeta).
\end{align*}
Now the Poisson integration for $H^1$ functions gives 
\begin{align*}
|f^{(1)}_\ell(\lambda)|^2\leqslant \int_{\mathbb T} |f^{(1)*}_\ell(\zeta)|^2 P_\lambda(\zeta) dm(\zeta)\leqslant  2 \mathrm{Re}\langle f^{(1)}_\ell \kappa^{\mu_1}_\lambda, f^{(1)}_\ell\rangle_{\mathcal D(\mu_1)}-\|f^{(1)}_\ell\|^2_{\mathcal D(\mu_1)}.
\end{align*}
Taking the sum over $\ell$ gives
\begin{align*}
\|L_{\lambda}(f)\|^2_{H^2(\mathbb D)}=\sum_{\ell=0}^\infty |f^{(1)}_\ell(\lambda)|^2 &\leqslant  2\sum_{\ell=0}^\infty  \mathrm{Re}\langle f^{(1)}_\ell \kappa^{\mu_1}_\lambda, f^{(1)}_\ell\rangle_{\mathcal D(\mu_1)}-\sum_{\ell=0}^\infty \|f^{(1)}_\ell\|^2_{\mathcal D(\mu_1)}. 
\end{align*}
Similarly, we can establish the inequality for $R_\lambda(f).$
\end{proof}
We are now ready to prove that every extremal function in $\mathcal D(\mu_1,\mu_2)$ is a multiplier (cf. \cite[Theorem~3.1]{RS1992} and \cite[Theorem~5.1(a)]{RS1994}).
\begin{proposition}\label{extremal-multiplier}
Every extremal function in $\mathcal D(\mu_1,\mu_2)$ is a multiplier of $\mathcal D(\mu_1,\mu_2).$ 
\end{proposition}
\begin{proof}
Let $\varphi$ be an extremal function in $\mathcal D(\mu_1,\mu_2)$. By Proposition~\ref{extremal-cyclic-subspace}, we have 
\begin{equation}\label{eqnextremal}
    \varphi \mathcal D(\nu_1,\nu_2)=[\varphi]\subseteq \mathcal D(\mu_1,\mu_2),
\end{equation}
with $d\nu_i(\zeta)=\rho_i(\zeta) d\mu_i(\zeta)$, $i=1,2,$ where
$$\rho_1(\zeta)=\sum_{\ell=0}^\infty |\varphi^{(1)*}_\ell(\zeta)|^2 ,~\rho_2(\zeta)=\sum_{k=0}^\infty |\varphi^{(2)*}_k(\zeta)|^2.$$ 
To prove that $\varphi$ is a multiplier of $\mathcal D(\mu_1,\mu_2)$, it suffices to show that  $\rho_i \in L^\infty(\mu_i)$ for $i=1,2.$ 
Indeed, if  $\rho_i \in L^\infty(\mu_i)$ for $i=1,2,$ then it is evident from the definition that $$\mathcal D(\mu_1,\mu_2)\subseteq \mathcal D(\nu_1,\nu_2)$$ 
(cf. \cite[Theorem~6.2]{R1991}). Hence, in view of \eqref{eqnextremal}, we get
\begin{align*}
\varphi \mathcal D(\mu_1,\mu_2)\subseteq \varphi \mathcal D(\nu_1,\nu_2)=[\varphi] \subseteq \mathcal D(\mu_1,\mu_2),
\end{align*}
establishing that $\varphi \in \mathrm{Mult}(\mathcal D(\mu_1,\mu_2)).$

 By applying Lemma~\ref{lemma:k-mu-mult}, for each $\lambda\in \mathbb D,$ we get  
  \begin{align}\label{eq:phi-lambda-bounded}
 \sum_{\ell=0}^\infty |\varphi^{(1)}_\ell(\lambda)|^2 &\leqslant  2\sum_{\ell=0}^\infty  \mathrm{Re}\langle \varphi^{(1)}_\ell \kappa^{\mu_1}_\lambda, \varphi^{(1)}_\ell \rangle_{\mathcal D(\mu_1)}-\sum_{\ell=0}^\infty \|\varphi^{(1)}_\ell\|^2_{\mathcal D(\mu_1)}. 
 \end{align}
 First, we show that the right hand side of the above inequality is finite. Since $\kappa^{\mu_1}_\lambda$  belongs to  $\mathrm{Mult}(\mathcal D(\mu_1),$ by Lemma~\ref{one-to-two-mult}, $\kappa^{\mu_1}_\lambda(z_1)$  belongs to  $\mathrm{Mult}(\mathcal D(\mu_1,\mu_2)).$ Thus $\kappa^{\mu_1}_\lambda(z_1) \varphi\in \mathcal D(\mu_1,\mu_2).$ Since $\varphi$ is extremal and $\kappa_{\lambda}^{\mu_1}(0)=1$, we have
 \allowdisplaybreaks
 \begin{eqnarray*}
&& \|\varphi\|_{\mathcal D(\mu_1,\mu_2)}^2\\
&=&\langle \kappa^{\mu_1}_\lambda(z_1)\varphi , \varphi \rangle_{\mathcal D(\mu_1,\mu_2)}+\langle \varphi , \kappa^{\mu_1}_\lambda(z_1)\varphi \rangle_{\mathcal D(\mu_1,\mu_2)}-\langle \varphi , \varphi \rangle_{\mathcal D(\mu_1,\mu_2)}\\
 &=& 2\mathrm{Re} \langle \kappa^{\mu_1}_\lambda(z_1)\varphi , \varphi \rangle_{\mathcal D(\mu_1,\mu_2)}-\langle \varphi , \varphi \rangle_{\mathcal D(\mu_1,\mu_2)}\\
 &=& 2\left(\mathrm{Re}\langle \kappa^{\mu_1}_\lambda(z_1)\varphi , \varphi \rangle_{H^2(\mathbb D^2)}+ \mathrm{Re} D^{(1)}_{\mu_1,\mu_2}(\kappa^{\mu_1}_\lambda(z_1)\varphi , \varphi)+\mathrm{Re} D^{(2)}_{\mu_1,\mu_2}(\kappa^{\mu_1}_\lambda(z_1)\varphi , \varphi)\right)\\
 &&-\left(\langle \varphi,\varphi\rangle_{H^2(\mathbb D^2)}+D^{(1)}_{\mu_1,\mu_2}(\varphi , \varphi)+D^{(2)}_{\mu_1,\mu_2}(\varphi , \varphi)\right)\\
 &=& 2\left(\sum_{\ell=0}^\infty \mathrm{Re} \langle \kappa^{\mu_1}_\lambda\varphi^{(1)}_\ell , \varphi^{(1)}_\ell \rangle_{H^2(\mathbb D)}+ \sum_{\ell=0}^\infty \mathrm{Re} D_{\mu_1}(\kappa^{\mu_1}_\lambda\varphi^{(1)}_\ell , \varphi^{(1)}_\ell)+ \mathrm{Re} D^{(2)}_{\mu_1,\mu_2}(\kappa^{\mu_1}_\lambda(z_1)\varphi , \varphi)\right)\\
 &&-\left(\sum_{\ell=0}^\infty \|\varphi^{(1)}_\ell\|^2_{H^2(\mathbb D)}+\sum_{\ell=0}^\infty D_{\mu_1}(\varphi^{(1)}_\ell)+D^{(2)}_{\mu_1,\mu_2}(\varphi , \varphi)\right) ~~(\text{by Lemma  \ref{weight-functions} and polarisation})\\
 &=& 2\sum_{\ell=0}^\infty \mathrm{Re} \langle \kappa^{\mu_1}_\lambda\varphi^{(1)}_\ell , \varphi^{(1)}_\ell \rangle_{\mathcal D(\mu_1)}-\sum_{\ell=0}^\infty \|\varphi^{(1)}_\ell\|^2_{\mathcal D(\mu_1)}+ 2\mathrm{Re} D^{(2)}_{\mu_1,\mu_2}(\kappa^{\mu_1}_\lambda(z_1)\varphi , \varphi)- D^{(2)}_{\mu_1,\mu_2}(\varphi , \varphi).
 \end{eqnarray*}
Note that 
 \begin{align*}
 &2\mathrm{Re} D^{(2)}_{\mu_1,\mu_2}(\kappa^{\mu_1}_\lambda(z_1)\varphi , \varphi)- D^{(2)}_{\mu_1,\mu_2}(\varphi , \varphi)\\
 &= \lim_{r\to 1^-} \int_{0}^{2\pi}\int_\mathbb D \left(2\mathrm{Re}\kappa^{\mu_1}_\lambda(re^{i\theta})-1\right)|\partial_2\varphi(re^{i\theta},z_2)|^2P_{\mu_2}(z_2)dA(z_2)\frac{d\theta}{2\pi} \overset{\eqref{2re-k-1}}{\geqslant} 0.
 \end{align*}
 This establishes
 \begin{align}\label{eq:kappa-mu12}
2\sum_{\ell=0}^\infty  \mathrm{Re}\langle \varphi^{(1)}_\ell \kappa^{\mu_1}_\lambda, \varphi^{(1)}_\ell \rangle_{\mathcal D(\mu_1)}-\sum_{\ell=0}^\infty \|\varphi^{(1)}_\ell\|^2_{\mathcal D(\mu_1)}\leqslant \|\varphi\|_{\mathcal D(\mu_1,\mu_2)}^2.
\end{align}
Combining \eqref{eq:phi-lambda-bounded} and \eqref{eq:kappa-mu12} yields 
\begin{align*}
\sum_{\ell=0}^\infty |\varphi^{(1)}_\ell(\lambda)|^2  \leqslant \|\varphi\|_{\mathcal D(\mu_1,\mu_2)}^2.
\end{align*} 
By applying Fatou's lemma we obtain
\begin{align*}
\sum_{\ell=0}^\infty |\varphi^{(1)*}_\ell(\zeta)|^2  \leqslant \|\varphi\|_{\mathcal D(\mu_1,\mu_2)}^2,~ \text{$\mu_1$-a.e. $\zeta\in \mathbb T.$}
\end{align*}
Similarly, one can show
\begin{align*}
\sum_{k=0}^\infty |\varphi^{(2)*}_k(\zeta)|^2  \leqslant \|\varphi\|_{\mathcal D(\mu_1,\mu_2)}^2,\quad \text{$\mu_2$-a.e. $\zeta\in \mathbb T.$}
\end{align*}
Hence  $\sum_{\ell=0}^\infty |\varphi^{(1)*}_\ell|^2$  and $\sum_{k=0}^\infty |\varphi^{(2)*}_k|^2$ belong to $ L^\infty(\mu_1)$ and $ L^\infty(\mu_2),$ respectively.
\end{proof}

\begin{lemma}\cite[Proposition~2.17]{BGV2026}\label{lem:BGV}
    Let $(T_1,T_2)$ be a left-inverse commuting tuple on $\mathcal H.$ If $\varphi\in \ker T_1^*\cap \ker T_2^*,$ then 
     \begin{equation*}
        \begin{split}
            \langle T_1^m\varphi, T_1^pT_2^{q+1}\varphi\rangle_{\mathcal D(\mu_1,\mu_2)}&=0,\quad m,p,q\geqslant 0,\\
            \langle T_2^n\varphi, T_1^{p+1}T_2^q\varphi\rangle_{\mathcal D(\mu_1,\mu_2)}&=0,\quad  n,p,q\geqslant 0.
        \end{split}
    \end{equation*}
\end{lemma}

We are now ready to prove Theorem~\ref{inv-sub-d-mu}. 
\begin{proof}[Proof of Theorem~\ref{inv-sub-d-mu}]
$\rm(i)$ Let $\mathcal S$ be a left-inverse commuting submodule of $\mathcal D(\mu_1,\mu_2)$. From \cite[Lemma~4.1]{BGV2025} it follows that the dimension of $\mathcal  S\ominus (z_1 \mathcal S+z_2 \mathcal S)$ is at least one. If possible, let there exist $\varphi,\psi \in \mathcal  S\ominus (z_1 \mathcal S+z_2 \mathcal S)$ such that $\varphi\perp \psi$ and $\|\varphi\|_{\mathcal D(\mu_1,\mu_2)}= \|\psi\|_{\mathcal D(\mu_1,\mu_2)}=1.$ Denote $$T_i=M_{z_i}|_\mathcal S,\quad i=1,2.$$  As $(T_1,T_2)$ are left-inverse commuting and $\varphi,\psi \in \ker T_1^*\cap \ker T_2^*,$  by Lemma~\ref{lem:BGV}, $\varphi$ and $\psi$ satisfy \eqref{extremal}, and hence are extremal functions in $\mathcal D(\mu_1,\mu_2).$ Applying Proposition~\ref{extremal-multiplier} we get that $\varphi$ and $\psi$ are multipliers of $\mathcal D(\mu_1,\mu_2).$ This implies that $\varphi \psi\in [\varphi]\cap [\psi].$ Since polynomials are dense in $\mathcal D(\mu_1,\mu_2)$ so there exist sequences $\{p_n\}_{n\geqslant 1}$ and $\{q_n\}_{n\geqslant 1}$ such that $p_n\rightarrow \psi$ and $q_n\rightarrow \varphi$ in $\mathcal D(\mu_1,\mu_2)$ as $n\to \infty.$ As both $\varphi$ and $\psi$ are multipliers so  $p_n\psi\rightarrow \varphi\psi$ and $q_n\varphi\rightarrow \varphi\psi$ in $\mathcal D(\mu_1,\mu_2)$ as $n\to \infty.$ Since $p_n$  and $q_n$ are polynomials so 
 \begin{align*}
 p_n\varphi=p_n(T_1,T_2)\varphi\underset{n\to\infty}{\longrightarrow} \varphi \psi, \quad \text{and}\quad  q_n\psi=q_n(T_1,T_2)\psi\underset{n\to\infty}{\longrightarrow} \varphi \psi.
 \end{align*}
For $i=1,2,$ denote $L_i=(T_i^*T_i)^{-1}T_i^*.$ Then $\ker T_i^*=\ker L_i$ and by assumption $T_iL_j=L_jT_i$ for $1\leqslant i\neq j\leqslant 2.$ Define 
\begin{align*}
P_i=I_\mathcal S-T_iL_i,\quad i=1,2.
\end{align*}
Then $P_i$ is a projection onto $\ker T_i^*,$ $i=1,2.$ Suppose 
$$p_n(z_1,z_2)=\sum_{i=0}^{a_n}\sum_{j=0}^{b_n}\widehat{p_n}(i,j)z_1^iz_2^j\text{  and  } q_n(z_1,z_2)=\sum_{i=0}^{c_n}\sum_{j=0}^{d_n}\widehat{q_n}(i,j)z_1^iz_2^j.$$
Fix $k\geqslant 0.$ As $L_1T_2=T_2L_1$ and $\varphi \in \ker L_1$ so  
\begin{align*}
L_1^kp_n(T_1,T_2)\varphi= \sum_{i=0}^{a_n}\sum_{j=0}^{b_n}\widehat{p_n}(i,j)L_i^kT_1^iT_2^j\varphi=\sum_{i=k}^{a_n}\sum_{j=0}^{b_n}\widehat{p_n}(i,j)T_1^{i-k}T_2^j\varphi.
\end{align*}
If $i>k$ then $$P_1T_1^{i-k}T_2^j\varphi= T_1^{i-k}T_2^j\varphi-T_1 T_1^{i-k-1}T_2^j\varphi=0,$$ and for $i=k,$ as $\ker T_1^*$ is $T_2$-reducing for left-inverse commuting tuple $(T_1,T_2)$ (see \cite[Corollary~4.2]{BGV2025}) so $P_1T_2^j\varphi=T_2^j\varphi.$ Hence
 \begin{align*}
P_1L_1^kp_n(T_1,T_2)\varphi=\sum_{i=k}^{a_n}\sum_{j=0}^{b_n}\widehat{p_n}(i,j)P_1T_1^{i-k}T_2^j\varphi= \sum_{j=0}^{b_n}\widehat{p_n}(k,j)T_2^j\varphi.
\end{align*}
Similarly, fixing $\ell\geqslant 0$ and applying $P_2L_2^\ell$ to the above equation gives 
\begin{align*}
P_2L_2^\ell P_1L_1^kp_n(T_1,T_2)\varphi= \sum_{j=0}^{b_n}\widehat{p_n}(k,j)P_2L_2^\ell T_2^j\varphi=\widehat{p_n}(k,\ell)\varphi.
\end{align*}
Hence for any pair $(k,\ell)\in \mathbb Z_+^2,$
\begin{equation}\label{l-k-l-in-phi}
P_2L_2^\ell P_1L_1^kp_n(T_1,T_2)\varphi=\widehat{p_n}(k,\ell)\varphi \underset{n\to\infty}{\longrightarrow} P_2L_2^\ell P_1L_1^k\varphi \psi.
\end{equation}
On the other hand, 
\begin{align*}
P_2L_2^\ell P_1L_1^kq_n(T_1,T_2)\psi=\widehat{q_n}(k,\ell)\psi,
\end{align*}
so 
\begin{equation}\label{l-k-l-in-psi}
P_2L_2^\ell P_1L_1^kq_n(T_1,T_2)\psi=\widehat{q_n}(k,\ell)\psi \underset{n\to\infty}{\longrightarrow} P_2L_2^\ell P_1L_1^k\varphi \psi.
\end{equation}
Since $\varphi\perp \psi$ so combining \eqref{l-k-l-in-phi} and \eqref{l-k-l-in-psi} gives $$P_2L_2^\ell P_1L_1^k\varphi \psi=0,\quad \text{for all } k,\ell\geqslant 0.$$
Let $h_k=P_1L_1^k\varphi \psi$ then $P_2L_2^\ell h_k=0$ for all $\ell\geqslant 0.$ Then for any integer $N\geqslant 1$ using the fact that $P_2=I_\mathcal S-T_2L_2$  we get
\begin{align*}
h_k-T_2^NL_2^Nh_k=\sum_{k=0}^{N-1} T_2^kL_2^kh_k-T_2^{k+1}L_2^{k+1}h_k=\sum_{k=0}^{N-1} T_2^k(I_\mathcal S-T_2L_2)L_2^kh_k=\sum_{k=0}^{N-1} T_2^kP_2L_2^kh_k=0
\end{align*} 
(see \cite[Eq.~2.2]{Sh2001-1}). Hence $h_k=T_2^NL_2^Nh_k$ for all $N\geqslant 0.$ Since $T_2$ is analytic and 
\begin{eqnarray*}
    h_k\in \cap_{N\geqslant 1} T_2^N\mathcal S.
\end{eqnarray*}This yields $0=h_k=P_1L_1^k\varphi \psi.$  Arguing similarly with $T_1,$ one can show that $\varphi\psi=0.$ Which implies one of $\varphi$ or $\psi$ is zero which is contradiction to the fact that $\|\varphi\|_{\mathcal D(\mu_1,\mu_2)}= \|\psi\|_{\mathcal D(\mu_1,\mu_2)}=1.$ So there cannot be two nonzero orthogonal vectors in $\mathcal  S\ominus (z_1 \mathcal S+z_2 \mathcal S)$ and hence the result.

$\rm(ii)$ $(\Rightarrow)$ 
Let $\mathcal S$ be a left-inverse commuting submodule of $\mathcal D(\mu_1,\mu_2).$ Then by part (i),
$$\dim \mathcal S\ominus (z_1\mathcal S+z_2\mathcal S)=1.$$
Let $T_i=M_{z_i}|_\mathcal S$ for $i=1,2$ and $\varphi \in \mathcal S\ominus (z_1\mathcal S+z_2\mathcal S)=\ker T_1^*\cap \ker T_2^*$ with $\|\varphi\|_{\mathcal D(\mu_1,\mu_2)}=1.$ By the assumption $(T_1,T_2)$ is left-inverse commuting, applying Lemma~\ref{lem:BGV} we get $\varphi$ is an extremal function in $\mathcal D(\mu_1,\mu_2)$.
Since $(T_1,T_2)$ is left-inverse commuting $2$-isometries, combining \cite[Theorem~1]{R1988} and  \cite[Theorem~4.3]{BGV2025} gives $\ker T^*_1\cap \ker T^*_2$ is a generating wandering subspace of $(T_1,T_2).$ Hence $$\mathcal S=\bigvee_{m,n\geqslant 0} T_1^mT_2^n(\ker T^*_1\cap \ker T^*_2)=[\varphi].$$
Finally, Proposition~\ref{extremal-cyclic-subspace} yields the result.

 $(\Leftarrow)$ Suppose there exists an extremal function $\varphi$ in $\mathcal D(\mu_1,\mu_2)$ such that $\mathcal S=\varphi\mathcal D(\nu_1,\nu_2)$ and $\nu_1,\nu_2$ are give by \eqref{nu-1-2-formula}. By applying Proposition~\ref{extremal-cyclic-subspace}, $\mathcal S=[\varphi]$ and $M_\varphi$ is an isometry from $\mathcal D(\nu_1,\nu_2)$ into $\mathcal D(\mu_1,\mu_2).$
Let $P$ denote the projection from $\mathcal D(\mu_1,\mu_2)$ onto $\varphi \mathcal D(\nu_1,\nu_2)$ then $P=M_\varphi M^*_\varphi.$  For $i=1,2,$ Let $T_i=M_{z_i}|_\mathcal S$ and $\tilde L_i$ be given by $$\tilde L_i=(T_i^*T_i)^{-1}T^*_i=(PM_{z_i}^*M_{z_i}P)^{-1}PM_{z_i}^*.$$
Next we show that $(T_1,T_2)$ is left-inverse commuting. For $1\leqslant i\neq j\leqslant 2,$
\begin{align}\label{l-tilde-i-t-j}
\tilde L_i T_j&=(PM_{z_i}^*M_{z_i}P)^{-1}PM_{z_i}^*M_{z_j}P =(M_\varphi M^*_\varphi M_{z_i}^*M_{z_i}M_\varphi M^*_\varphi)^{-1}M_\varphi M^*_\varphi M_{z_i}^*M_{z_j}M_\varphi M^*_\varphi \notag \\
&=(M_\varphi  M_{z_i}^* M^*_\varphi M_\varphi M_{z_i} M^*_\varphi)^{-1}M_\varphi M_{z_i}^* M^*_\varphi M_\varphi M_{z_j} M^*_\varphi =(M_\varphi  M_{z_i}^* M_{z_i} M^*_\varphi)^{-1}M_\varphi (M_{z_i}^*M_{z_j}) M^*_\varphi.
\end{align}
Since $(M_{z_1},M_{z_2})$ is left-inverse commuting so for $1\leqslant i\neq j\leqslant 2,$
$$M_{z_i}^*M_{z_j}= M_{z_i}^*M_{z_i}M_{z_j}(M_{z_i}^*M_{z_i})^{-1}M_{z_i}^*.$$
Now substituting this in \eqref{l-tilde-i-t-j} gives 
\begin{eqnarray}\label{l-i-t-j}
\tilde L_i T_j&=&(M_\varphi  M_{z_i}^* M_{z_i} M^*_\varphi)^{-1} M_\varphi M_{z_i}^* M_{z_i} M_{z_j}(M_{z_i}^*M_{z_i})^{-1} M_{z_i}^* M^*_\varphi \notag\\
&=& (M_\varphi  M_{z_i}^* M_{z_i} M^*_\varphi)^{-1}(M_\varphi M_{z_i}^* M_{z_i}M^*_\varphi) M_\varphi M_{z_j}(M_{z_i}^*M_{z_i})^{-1}M_{z_i}^* M^*_\varphi \notag \\
&=& M_\varphi M_{z_j}(M_{z_i}^*M_{z_i})^{-1}M_{z_i}^* M^*_\varphi.
\end{eqnarray}
On the other hand, for $1\leqslant i\neq j\leqslant 2,$\allowdisplaybreaks
\begin{eqnarray}\label{t-j-l-i}
T_j \tilde L_i &=& M_{z_j}P(PM_{z_i}^*M_{z_i}P)^{-1}PM_{z_i}^* \notag \\
&=& M_{z_j}M_\varphi M^*_\varphi (M_\varphi M^*_\varphi M_{z_i}^*M_{z_i}M_\varphi M^*_\varphi )^{-1}M_\varphi M^*_\varphi M_{z_i}^* \notag \\
&=& M_{z_j}M_\varphi M^*_\varphi (M_\varphi  M_{z_i}^* M^*_\varphi M_\varphi M_{z_i} M^*_\varphi )^{-1}M_\varphi M^*_\varphi M_{z_i}^* \notag \\
&=& M_{z_j}M_\varphi (M_{z_i}^* M_{z_i})^{-1} M_{z_i}^* M_{z_i} M^*_\varphi (M_\varphi  M_{z_i}^* M_{z_i} M^*_\varphi )^{-1}M_\varphi M^*_\varphi M_{z_i}^* \notag \\
&=& M_{z_j}M_\varphi (M_{z_i}^* M_{z_i})^{-1}  M^*_\varphi ( M_\varphi M_{z_i}^* M_{z_i} M^*_\varphi) (M_\varphi  M_{z_i}^* M_{z_i} M^*_\varphi )^{-1}M_\varphi M^*_\varphi M_{z_i}^* \notag \\
&=& M_{z_j}M_\varphi (M_{z_i}^* M_{z_i})^{-1}  M^*_\varphi M_\varphi M^*_\varphi M_{z_i}^* \notag \\
&=& M_\varphi M_{z_j}(M_{z_i}^* M_{z_i})^{-1} M_{z_i}^* M^*_\varphi. 
\end{eqnarray} 
Thus \eqref{l-i-t-j} and \eqref{t-j-l-i} shows that $\mathcal S$ is a left-inverse commuting submodule.
\end{proof}

Two submodules $\mathcal S_1$ and $\mathcal S_2$ of $\mathcal D(\mu_1,\mu_2)$ are said to be unitarily equivalent  if the commuting tuple $(M_{z_1}|_{\mathcal S_1},M_{z_2}|_{\mathcal S_1})$ is unitarily equivalent to $(M_{z_1}|_{\mathcal S_2},M_{z_2}|_{\mathcal S_2}).$
\begin{theorem}\label{thm:uni-equiv}
Let $\mathcal S_1$ and $\mathcal S_2$ be nonzero left-inverse commuting
submodules of $\mathcal D(\mu_1,\mu_2)$. For $i=1,2$, suppose $\mathcal S_i=\varphi_i\mathcal D(\nu_1^{(i)},\nu_2^{(i)})$, where  $\varphi_i$
is a unit vector in $\mathcal S_i\ominus(z_1\mathcal S_i+z_2\mathcal S_i)$ $($as given in Theorem~\ref{inv-sub-d-mu}$\rm{(ii)}$$)$. Then $\mathcal S_1$ and
$\mathcal S_2$ are unitarily equivalent  if and only if
$$\nu_1^{(1)}=\nu_1^{(2)}\quad\text{and}\quad\nu_2^{(1)}=\nu_2^{(2)}.$$
\end{theorem}

\begin{proof}
By Theorem~\ref{inv-sub-d-mu}(ii), for $i=1,2,$ the following map
$$J_i: \mathcal D(\nu_1^{(i)},\nu_2^{(i)})\rightarrow\mathcal S_i,\qquad J_i f=\varphi_i f,~~f\in \mathcal D(\nu_1^{(i)},\nu_2^{(i)})$$ is an unitary. If $W:\mathcal S_1\rightarrow\mathcal S_2$ is a unitary module map, then $$J_2^{-1}WJ_1:\mathcal D(\nu_1^{(1)},\nu_2^{(1)})\rightarrow \mathcal D(\nu_1^{(2)},\nu_2^{(2)})$$ is also a unitary module map. Hence \cite[Theorem~6.4]{BCG2025} gives $\nu_j^{(1)}=\nu_j^{(2)}$ for $j=1,2$.

Conversely, if the two pairs of measures agree, then
$\mathcal D(\nu_1^{(1)},\nu_2^{(1)})=\mathcal D(\nu_1^{(2)},\nu_2^{(2)})$. Thus $J_2J_1^{-1}$ is a unitary module
map from $\mathcal S_1$ onto $\mathcal S_2$, given by
$\varphi_1 f\mapsto\varphi_2 f$.
\end{proof}

\section{Extremal functions of \texorpdfstring{$\mathcal D(\mathbb D^2)$}{}}\label{section-5}
In case of $\mathcal D(\mathbb D^2),$ we show that a function is extremal if and only if it is a constant multiple  of some monomial. We start with the observation that extremal functions of $\mathcal D(\mathbb D^2)$ are constant multiple of inner functions.
\begin{lemma}\label{phi-inner}
Let $\varphi\in \mathcal D(\mathbb D^2)$ be an extremal function. Then $\varphi$ is a constant multiple of inner function. 
\end{lemma}
\begin{proof}
For  $f,g\in \mathcal D(\mathbb D^2)$, we have
\begin{align*}
\langle f,g\rangle_{\mathcal D(\mathbb D^2)}=\sum_{k,\ell=0}^\infty (k+\ell+1) \hat{f}(k,\ell)\overline{\hat{g}(k,\ell)}.
\end{align*}
Thus 
\begin{equation}\label{richter-formula}
\langle z_if,z_ig\rangle_{\mathcal D(\mathbb D^2)}-\langle f,g\rangle_{\mathcal D(\mathbb D^2)}=\langle f,g\rangle_{H^2(\mathbb D^2)},~i=1,2.
\end{equation}
As $\varphi\in \mathcal D(\mathbb D^2)\subseteq H^2(\mathbb D^2)$,  $\varphi^*(\zeta_1,\zeta_2)$ exists a.e. on $\mathbb T^2$ and $\varphi^*\in L^2(\mathbb T^2)$ (see \cite[Theorem~3.3.5]{Ru1969}). Now consider the function 
\begin{align*}
u(\zeta_1,\zeta_2):=|\varphi^*(\zeta_1,\zeta_2)|^2, \quad  (\zeta_1,\zeta_2)\in \mathbb T^2.
\end{align*}
Then for $k,\ell\geqslant 0,$ 
\begin{align*}
\widehat{u}(k,\ell)&=\int_{\mathbb T^2} u(\zeta_1,\zeta_2)\overline{\zeta}_1^k\overline{\zeta}_2^\ell dm(\zeta_1) dm(\zeta_2)\\
&=\int_{\mathbb T^2} |\varphi^*(\zeta_1,\zeta_2)|^2\overline{\zeta}_1^k\overline{\zeta}_2^\ell dm(\zeta_1) dm(\zeta_2)\\
&= \int_{\mathbb T^2} \varphi^*(\zeta_1,\zeta_2)\overline{\zeta_1^k\zeta_2^\ell \varphi^*(\zeta_1,\zeta_2)}dm(\zeta_1) dm(\zeta_2)= \langle \varphi, z_1^kz_2^\ell\varphi\rangle_{H^2(\mathbb D^2)}.
\end{align*}
Similarly, one can show that for $k,\ell\geqslant 0,$
\begin{align*}
\widehat{u}(k,-\ell)=\langle z_2^\ell\varphi, z_1^k\varphi\rangle_{H^2(\mathbb D^2)}.
\end{align*}
Let $k\geqslant 0$ and $\ell>0.$  Since $\varphi$ is extremal, using \eqref{richter-formula} and \eqref{extremal}  we get
\begin{align*}
\widehat{u}(k,\ell)&=\langle \varphi, z_1^k z_2^\ell\varphi\rangle_{H^2(\mathbb D^2)}=\langle z_1\varphi,z_1^{k+1}z_2^\ell \varphi\rangle_{\mathcal D(\mathbb D^2)}-\langle \varphi,z_1^{k}z_2^\ell\varphi\rangle_{\mathcal D(\mathbb D^2)}=0,\\
\widehat{u}(k,-\ell)&=\langle z_2^\ell\varphi, z_1^k\varphi\rangle_{H^2(\mathbb D^2)}=\langle z_1z_2^\ell \varphi,z_1^{k+1}\varphi\rangle_{\mathcal D(\mathbb D^2)}-\langle z_2^\ell\varphi,z_1^{k}\varphi\rangle_{\mathcal D(\mathbb D^2)}=0,\\
\widehat{u}(k+1,0)&=\langle \varphi, z_1^{k+1}\varphi\rangle_{H^2(\mathbb D^2)}=\langle z_2\varphi,z_1^{k+1}z_2\varphi\rangle_{\mathcal D(\mathbb D^2)}-\langle \varphi,z_1^{k+1}\varphi\rangle_{\mathcal D(\mathbb D^2)}=0.
\end{align*}
Thus $\widehat{u}(k,\ell)=0$ for all $(k,\ell)\in (\mathbb Z_+\times \mathbb Z)\backslash \{(0,0)\}.$ Since $u$ is real valued so $$\widehat{u}(-k,-\ell)=\overline{\widehat{u}(k,\ell)},\quad (k,\ell)\in \mathbb Z\times \mathbb Z.$$
Hence $\widehat{u}(k,\ell)=0$ for $(k,\ell)\in \mathbb Z^2\backslash \{(0,0)\}$. Hence $u$ is constant a.e. on $\mathbb T^2$ and 
\begin{equation*}
u(\zeta_1,\zeta_2)=|\varphi^*(\zeta_1,\zeta_2)|^2=\widehat{u}(0,0)=\langle \varphi, \varphi\rangle_{H^2(\mathbb D^2)}=\|\varphi \|^2_{H^2(\mathbb D^2)}.
\end{equation*}
Since $\varphi\in H^2(\mathbb D^2)$,  it follows that
\begin{align*}
|\varphi(\lambda_1,\lambda_2)|^2\leqslant \int_{\mathbb T^2} P_{\lambda_1}(\zeta_1)P_{\lambda_2}(\zeta_2)|\varphi^*(\zeta_1,\zeta_2)|^2dm(\zeta_1) dm(\zeta_2)= \|\varphi \|^2_{H^2(\mathbb D^2)}.
\end{align*}
Thus $\varphi\in H^\infty(\mathbb D^2).$ Hence $\varphi$ is a constant multiple of an inner function.
\end{proof}

\begin{lemma}\label{partial-f-in-h2}
    Let $f\in \mathcal D(\mathbb D^2).$ Then for each $\lambda=(\lambda_1,\lambda_2) \in \mathbb D^2,$ the functions 
$$w\mapsto\partial_1 f(\lambda_1,w)\text{ and  } w\mapsto \partial_2f(w,\lambda_2)\text{ belongs to } H^2(\mathbb D).$$
\end{lemma}
\begin{proof}
 Let $f(z_1,z_2)=\sum_{k=0}^\infty \sum_{\ell=0}^\infty \hat{f}(k,\ell)z_1^kz_2^\ell.$ Then 
 \begin{align*}
 D_{m,m}(f)=\sum_{k=0}^\infty \sum_{\ell=0}^\infty (k+\ell)|\hat{f}(k,\ell)|^2<\infty.
 \end{align*}
Note that  
 $$\partial_2f(w,\lambda_2)=\sum_{k=0}^\infty \left( \sum_{\ell=1}^\infty \ell\hat{f}(k,\ell)\lambda_2^{\ell-1}\right)w^k.$$
 Observe that 
 \begin{align*}
  \sum_{k=0}^\infty \left |\sum_{\ell=1}^\infty \ell\hat{f}(k,\ell)\lambda_2^{\ell-1}\right|^2
 &= \sum_{k=0}^\infty \left|\sum_{\ell=1}^\infty \left(\sqrt{\ell}\lambda_2^{\ell-1}\right)\left(\sqrt{\ell}\hat{f}(k,\ell)\right)\right|^2\\
 &\leqslant  \sum_{k=0}^\infty \left(\sum_{\ell=1}^\infty \ell|\lambda_2|^{2(\ell-1)}\right)\left(\sum_{\ell=1}^\infty \ell|\hat{f}(k,\ell)|^2\right)\\
 &\leqslant  \frac{1}{(1-|\lambda_2|^2)^2}  D_{m,m}(f)<\infty.
 \end{align*}
 Hence $w\mapsto \partial_2f(w,\lambda_2)\in H^2(\mathbb D).$ Similarly, one can show that $w\mapsto\partial_1 f(\lambda_1,w)\in H^2(\mathbb D).$
\end{proof}

\begin{lemma}\label{euler-lemma}
   Let $f$ be a nonconstant holomorphic function on $\mathbb D^2.$ Suppose there exists a nonnegative integer $N$ such that $$z_1\partial_1f(z_1,z_2)+z_2\partial_2f(z_1,z_2)=Nf(z_1,z_2),\quad (z_1,z_2)\in \mathbb D^2.$$ Then $f$ is a homogeneous polynomial of degree $N.$
\end{lemma}
\begin{proof}
Suppose $f(z_1,z_2)=\sum_{k=0}^\infty \sum_{\ell=0}^\infty \hat{f}(k,\ell)z_1^k z_2^\ell.$ Then
\begin{align*}
   (z_1\partial_1+z_2\partial_2)f(z_1,z_2)= \sum_{k=0}^\infty \sum_{\ell=0}^\infty (k+\ell)\hat{f}(k,\ell)z_1^k z_2^\ell.
\end{align*}
Since $z_1\partial_1f+z_2\partial_2f=Nf,$ comparing coefficients on the both sides gives  
$$(k+\ell-N)\hat{f}(k,\ell)=0,\quad \text{for all } k,\ell\geqslant 0.$$
Hence  $\hat{f}(k,\ell)=0$ whenever $k+\ell\neq N.$ Therefore $f$ becomes
$$f(z_1,z_2)=\sum_{k=0}^{N}\hat{f}(k,N-k) z_1^k z_2^{N-k}.$$
This completes the proof.
\end{proof}

\begin{lemma}\label{lem:f-in-dmm}
    Let $f\in \mathcal D(\mathbb D^2)$. Then there exists $E\subset \mathbb T$ such that $m(E)=0$ and
    \begin{eqnarray}\label{B-d-1-integral}
D^{(1)}_{m,m}(f)= \int_{\mathbb T\backslash E} D_m(f(\cdot,\zeta))dm(\zeta),~~
D^{(2)}_{m,m}(f)= \int_{\mathbb T\backslash E} D_m(f(\zeta, \cdot))dm(\zeta).
\end{eqnarray}
\end{lemma}
\begin{proof}
    Considering $H^2(\mathbb D^2)$ as $H^2(\mathbb D)$-valued Hardy space $H^2(\mathbb D; H^2(\mathbb D))$ on the unit disc, and applying \cite[Theorem~4.6.A, p.~84]{RR1997} we get that there exists $E\subset \mathbb T$ such that $m(E)=0$ and  
    $$f(\zeta,\cdot), f(\cdot,\zeta)\in H^2(\mathbb D),\quad \zeta\in\mathbb T\backslash E.$$ 
    Note that, for each $z_2\in \mathbb D,$ $r\mapsto \int_\mathbb T|\partial_2 f(r\zeta,z_2)|^2 dm(\zeta)$ is an increasing function (see \cite[Remark 17.8 (a)]{Ru1969}). Hence, by applying Fubini-Tonelli theorem and the monotone convergence theorem, we get
\begin{align*}
D^{(2)}_{m,m}(f)&=\lim_{r\to 1^-}\int_\mathbb T \int_\mathbb D |\partial_2 f(r\zeta,z_2)|^2 dA(z_2)dm(\zeta)\\
&= \lim_{r\to 1^-}\int_\mathbb D\left( \int_\mathbb T|\partial_2 f(r\zeta,z_2)|^2 dm(\zeta)\right)dA(z_2)\qquad(\text{by Fubini-Tonelli})\\
&= \int_\mathbb D\left(\lim_{r\to 1^-} \int_\mathbb T|\partial_2 f(r\zeta,z_2)|^2 dm(\zeta)\right)dA(z_2).\qquad(\text{by monotone convergence})
\end{align*}
By Lemma~\ref{partial-f-in-h2}, we know that for each $z_2\in \mathbb D,$ $w\mapsto\partial_2f(w,z_2)\in H^2(\mathbb D)$ so 
$$\lim_{r\to 1^-} \int_\mathbb T|\partial_2 f(r\zeta,z_2)|^2 dm(\zeta)=\int_\mathbb T|\partial_2 f(\zeta,z_2)|^2 dm(\zeta).$$
Hence again applying Fubini-Tonelli theorem gives
\begin{eqnarray}\label{B-d-2-integral-new}
D^{(2)}_{m,m}(f)= \int_\mathbb D\left(\int_{\mathbb T}|\partial_2 f(\zeta,z_2)|^2 dm(\zeta)\right)dA(z_2)= \int_{\mathbb T}\left(\int_\mathbb D|\partial_2 f(\zeta,z_2)|^2 dA(z_2)\right)dm(\zeta).
\end{eqnarray}
As for each $\zeta\in \mathbb T\backslash E,$ $f(\zeta,\cdot)$ belongs to $H^2(\mathbb D)$ so \eqref{B-d-2-integral-new} takes the form 
\begin{eqnarray*}
D^{(2)}_{m,m}(f)= \int_{\mathbb T\backslash E} D_m(f(\zeta,\cdot))dm(\zeta).
\end{eqnarray*}
Similarly, one can derive that
\begin{eqnarray*}
D^{(1)}_{m,m}(f)= \int_{\mathbb T\backslash E} D_m(f(\cdot,\zeta))dm(\zeta).
\end{eqnarray*}
Hence the result.
\end{proof}
Carleson \cite[Theorem]{C1960} characterized the inner functions in the classical Dirichlet space $\mathcal D(\mathbb D)$ as the finite Blaschke products. On the bidisk, every rational inner function belongs to $\mathcal D(\mathbb D^2)$ (see \cite[Corollary~2.4]{BS2026}). Below we prove the converse.
\begin{proposition}\label{prop:RIF}
    Only inner functions in $\mathcal D(\mathbb D^2)$ are rational inner functions.
\end{proposition}
\begin{proof}
    From \cite[Corollary~2.4]{BS2026}, we know that every rational inner function belongs $\mathcal D(\mathbb D^2)$. Next we show that these are all inner function in $\mathcal D(\mathbb D^2)$. Let $f$ be an inner function in $\mathcal D(\mathbb D^2)$. By Lemma~\ref{lem:f-in-dmm}, there exists $E\subset \mathbb T$ such that $m(E)=0$ and  
    \begin{align}\label{eq:dmmphi}
        D_{m,m}(f)= \int_{\mathbb T\backslash E} D_m(f(\cdot,\zeta))dm(\zeta)+ \int_{\mathbb T\backslash E} D_m(f(\zeta, \cdot))dm(\zeta)<\infty.
    \end{align}
    This implies that there exists $F\subset\mathbb T$ containing $E$ of Lebesgue measure zero such that for $\zeta \in \mathbb T\backslash F,$ $f(\zeta,\cdot)$ and $f(\cdot,\zeta)$ are inner functions on $\mathbb D$ and have finite Dirichlet integrals. From \cite[Theorem]{C1960}, we know that the only inner functions on $\mathbb D$ with finite Dirichlet integrals are finite Blaschke products. So for $\zeta \in \mathbb T\backslash F,$ $f(\zeta,\cdot)$ and $f(\cdot,\zeta)$ are finite Blaschke products. From \cite[Theorem~4]{Ah1971}, we know that 
    \begin{align*}
        &\text{if a two-variable inner function $\phi$ has an analytic extension across $\{\zeta\}\times \mathbb T$ for some $\zeta\in \mathbb T$} \\ &\text{and across $\mathbb T\times \{\eta\}$ for some $\eta\in \mathbb T,$ then $\phi$ is a rational inner function.}
    \end{align*}  
    Since  $f(\zeta,\cdot)$ and $f(\cdot,\zeta)$ are finite Blaschke products for each $\zeta \in \mathbb T\backslash F,$  by the above result, $f$ is a rational inner function.
\end{proof}

\begin{lemma}\label{lem:d-mm-pq}
    Let $f$ be an inner function on $\mathbb D^2.$  Suppose $F\subset \mathbb T$ with $m(F)=0$ and for each $\zeta_1,\zeta_2\in \mathbb T\backslash F,$  the functions $f(\cdot,\zeta_2)$ and $f(\zeta_1,\cdot)$  are holomorphic in a neighbourhood of $\overline{\mathbb D}.$ Then for any polynomials $p$ and $q,$ 
\begin{equation*}
\langle pf,qf\rangle_{\mathcal D(\mathbb D^2)}= \langle p,q \rangle_{\mathcal D(\mathbb D^2)}+ \int_{\mathbb T\backslash F\times \mathbb T\backslash F } p\overline{q}h(\zeta_1,\zeta_2)dm(\zeta_2)dm(\zeta_1),
\end{equation*}
where $h(\zeta_1,\zeta_2)=(\zeta_1 \partial_1 f + \zeta_2 \partial_2 f)(\zeta_1,\zeta_2)\overline{f(\zeta_1,\zeta_2)},$ $\zeta_1,\zeta_2\in \mathbb T\backslash F.$
\end{lemma}
\begin{proof}
     Let $g,h$ be two functions holomorphic in an neighbourhood of $\overline{\mathbb D}$ then 
\begin{equation}\label{one-variable-diri-int}
D_m(g,h)=\sum_{k=1}^\infty k\hat{g}(k)\overline{\hat{h}(k)}=\langle zg',h\rangle_{H^2(\mathbb D)}=\int_\mathbb T \eta g'(\eta)\overline{h(\eta)} dm(\eta).
\end{equation}
Since for each $\zeta_1,\zeta_2\in \mathbb T\backslash F,$  we know that $f(\cdot,\zeta_2)$ and $f(\zeta_1,\cdot)$  are finite Blaschke products, so they are holomorphic in a neighbourhood of $\overline{\mathbb D}.$
Thus for polynomials $p,q,$ 
\begin{eqnarray*}
&&D_{m,m}(pf,qf)\\
&=&D^{(1)}_{m,m}(pf,qf)+D^{(2)}_{m,m}(pf,qf)\\
&=& \int_{\mathbb T\backslash F} D_m((pf)(\cdot,\zeta_2),(qf)(\cdot,\zeta_2))dm(\zeta_2)+\int_{\mathbb T\backslash F} D_m((pf)(\zeta_1,\cdot),(qf)(\zeta_1,\cdot))dm(\zeta_1)\\
&\overset{\eqref{one-variable-diri-int}}{=}& \int_{\mathbb T\backslash F} \int_{\mathbb T} \zeta_1 \partial_1 (pf)(\zeta_1,\zeta_2) \overline{(qf)(\zeta_1,\zeta_2)}dm(\zeta_1)dm(\zeta_2)\\
&&+\int_{\mathbb T\backslash F} \int_{\mathbb T} \zeta_2 \partial_2 (pf)(\zeta_1,\zeta_2) \overline{(qf)(\zeta_1,\zeta_2)}dm(\zeta_2)dm(\zeta_1)\\
&=& \int_{\mathbb T\backslash F\times \mathbb T\backslash F }(\zeta_1 \partial_1 + \zeta_2 \partial_2)(pf)(\zeta_1,\zeta_2)\overline{(qf)(\zeta_1,\zeta_2)}dm(\zeta_2)dm(\zeta_1).
\end{eqnarray*}
Note that 
\begin{align*}
(z_1 \partial_1 + z_2 \partial_2)(pf)= f(z_1 \partial_1 p + z_2 \partial_2 p)+p(z_1 \partial_1 f + z_2 \partial_2 f).
\end{align*}
Since $f$ is inner so 
\allowdisplaybreaks
\begin{eqnarray*}
D_{m,m}(pf,qf)
&=& \int_{\mathbb T\backslash F\times \mathbb T\backslash F}(\zeta_1 \partial_1 p+ \zeta_2 \partial_2p)(\zeta_1,\zeta_2)\overline{q(\zeta_1,\zeta_2)}dm(\zeta_2)dm(\zeta_1)\\
&&+\int_{\mathbb T\backslash F\times \mathbb T\backslash F} p(\zeta_1,\zeta_2)\overline{q(\zeta_1,\zeta_2)}(\zeta_1 \partial_1 f + \zeta_2 \partial_2f)(\zeta_1,\zeta_2)\overline{f(\zeta_1,\zeta_2)}dm(\zeta_2)dm(\zeta_1)\\
&=& \int_{\mathbb T\backslash F\times \mathbb T\backslash F}\zeta_1 \partial_1 p(\zeta_1,\zeta_2)\overline{q(\zeta_1,\zeta_2)}dm(\zeta_2)dm(\zeta_1)\\
&&+\int_{\mathbb T\backslash F\times \mathbb T\backslash F} \zeta_2 \partial_2p(\zeta_1,\zeta_2)\overline{q(\zeta_1,\zeta_2)}dm(\zeta_2)dm(\zeta_1)\\
&& + \int_{\mathbb T\backslash F\times \mathbb T\backslash F} p\overline{q}hdm(\zeta_2)dm(\zeta_1)\\
&\overset{\eqref{one-variable-diri-int}}{=}& \int_{\mathbb T\backslash F} D_m(p(\cdot,\zeta_2),q(\cdot,\zeta_2))dm(\zeta_2)+\int_{\mathbb T\backslash F} D_m(p(\zeta_1,\cdot),q(\zeta_1,\cdot))dm(\zeta_1)\\
&& + \int_{\mathbb T\backslash F\times \mathbb T\backslash F} p\overline{q}hdm(\zeta_2)dm(\zeta_1)\\
&\overset{\eqref{B-d-1-integral}}{=}& D_{m,m}(p,q)+ \int_{\mathbb T\backslash F\times \mathbb T\backslash F} p\overline{q}hdm(\zeta_2)dm(\zeta_1).
\end{eqnarray*}
As $f$ is inner so $\langle pf,qf\rangle_{H^2(\mathbb D^2)}=\langle p,q\rangle_{H^2(\mathbb D^2)}$ and hence 
\begin{equation*}
\langle pf,qf\rangle_{\mathcal D(\mathbb D^2)}= \langle p,q \rangle_{\mathcal D(\mathbb D^2)}+ \int_{\mathbb T\backslash F\times \mathbb T\backslash F } p\overline{q}h(\zeta_1,\zeta_2)dm(\zeta_2)dm(\zeta_1).
\end{equation*}
\end{proof}

\begin{lemma}\label{lem:partial-f-bounded}
    Let $f\in H^\infty(\mathbb D^2).$ Suppose $F\subset \mathbb T$ with $m(F)=0.$ If for every $\zeta_1,\zeta_2\in \mathbb T\backslash F,$ 
    \begin{itemize}
        \item[\rm(i)] the functions $f(\cdot,\zeta_2)$ and $f(\zeta_1,\cdot)$ are holomorphic in a neighbourhood of $\overline{\mathbb D},$
        \item[\rm(ii)] there exists a constant $c\geqslant0$ such that  
      $$|\partial_i f(\zeta_1,\zeta_2)|\leq c,~~i=1,2,$$ 
    \end{itemize} then $\partial_i f \in H^\infty(\mathbb D^2)$ with 
          $$|\partial_i f(z_1,z_2)|\leq c,~~(z_1,z_2)\in \mathbb D^2,~~i=1,2.$$
\end{lemma}
\begin{proof}
Fix $\zeta_2 \in \mathbb T\backslash F.$
Consider the function $$b_{\zeta_2}(w)=f(w,\zeta_2),\quad w\in \mathbb D.$$ Then $b_{\zeta_2}$  is holomorphic in a neighbourhood of $\overline{\mathbb D}$. Hence $b'_{\zeta_2}$ is also holomorphic in a neighbourhood of $\overline{\mathbb D}.$ We know that $|b'_{\zeta_2}(\zeta_1)|\leq c$ for $\zeta_1\in \mathbb T\backslash F.$ Since $b'_{\zeta_2}$ is  holomorphic in a neighbourhood of $\overline{\mathbb D},$  $|b'_{\zeta_2}(\zeta_1)|\leq c,$ for every $\zeta_1\in \mathbb T.$
By maximum modulus principle we have
\begin{equation}\label{eq:b-zeta}
    |b'_{\zeta_2}(z_1)|\leq c,\quad z_1\in \mathbb D.
\end{equation}
Fix $z_1\in \mathbb D.$ Consider the function $$g_{z_1}(w)=\partial_1 f(z_1,w),\quad w\in \mathbb D.$$ For each $w\in \mathbb D,$ using Cauchy integral formula we have
\begin{align*}
g_{z_1}(w)=\partial_1 f(z_1,w)=\frac{1}{2\pi i}\int_{|z-z_1|=\rho} \frac{f(z,w)}{(z-z_1)^2}dz,~~\text{where }\rho=\frac{1-|z_1|}{2}.
\end{align*}
Since $f\in H^\infty(\mathbb D),$ 
\begin{align*}
|g_{z_1}(w)|=|\partial_1 f(z_1,w)|\leqslant \frac{1}{2\pi}2\pi \rho \frac{\|f\|_{\infty,\mathbb D^2}}{\rho^2}\leqslant \frac{1}{\rho}= \frac{2 \|f\|_{\infty,\mathbb D^2}}{1-|z_1|}.
\end{align*}
This shows that $g_{z_1}\in H^\infty (\mathbb D).$ Now 
\begin{align*}
    g_{z_1}^*(\zeta_2)=\lim_{r\to 1^-} g_{z_1}(r\zeta_2)= \lim_{r\to 1^-} \partial_1 f(z_1,r\zeta_2)=b'_{\zeta_2}(z_1).
\end{align*}
The last equality is due to local uniform convergence of holomorphic function. 
By \eqref{eq:b-zeta},  $$|g_{z_1}^*(\zeta_2)|\leq d.$$  This is true for each $\zeta_2\in \mathbb T\backslash F.$ As $g_{z_1}\in H^\infty (\mathbb D),$ using the Poisson integral formula gives
\begin{align*}
|\partial_1 f(z_1,w)|=|g_{z_1}(w)|\leqslant \int_\mathbb T P_{z_2}(\zeta_2)|g_{z_1}^*(\zeta_2)| dm(\zeta_2)\leqslant d,\quad w\in \mathbb D.
\end{align*}
This is true for each $z_1\in \mathbb D.$ Hence $\partial_1 f\in H^\infty (\mathbb D^2).$ Similarly, $\partial_2 f\in H^\infty (\mathbb D^2).$ 
\end{proof}
\begin{proof}[Proof of Theorem~\ref{extremals-are-monomials}]
Since monomials are orthogonal in $\mathcal D(\mathbb D^2)$, it follows that every monomial (upon normalizing) is an extremal function.

For the converse, suppose $\varphi\in\mathcal D(\mathbb D^2)$ is extremal. By Lemma~\ref{phi-inner}, we know that $\varphi$ is a constant multiple of an inner function with  $|\varphi^*(\zeta_1,\zeta_2)|=\|\varphi\|_{H^2(\mathbb D^2)}$  a.e. on $\mathbb T^2.$ Let us define $$\Phi(z_1,z_2)=\frac{\varphi(z_1,z_2)}{\|\varphi\|_{H^2(\mathbb D^2)}},\quad (z_1,z_2)\in \mathbb D^2.$$  Then $\Phi$ is an inner function on $\mathbb D^2.$ 
We first show that there exists a nonnegative integer $N$ such that 
\begin{equation}\label{eqnextremalD}
z_1\partial_1\Phi(z_1,z_2)+z_2\partial_2\Phi(z_1,z_2)=N\Phi(z_1,z_2),~(z_1,z_2)\in\mathbb D^2.
\end{equation}
Assume that \eqref{eqnextremalD} holds. Then  by Lemma~\ref{euler-lemma}, $\Phi$ becomes a homogeneous polynomial of degree $N.$ It is known that polynomial inner functions are unimodular constant times monomials (see \cite[Lemma~3.1]{DS2026}). Since $\Phi$ is inner so there exists  $\theta\in \mathbb R$  and $0\leqslant k_0\leqslant N$ such that $$\Phi(z_1,z_2)=e^{i\theta}z_1^{k_0}z_2^{N-k_0}.$$
By definition $\varphi=\|\varphi\|_{H^2(\mathbb D)}\Phi$ and $\|\varphi\|_{\mathcal D(\mathbb D^2)}=1.$ This implies $$1=\|\varphi\|^2_{H^2(\mathbb D)}(1+N).$$ Hence we get $$\varphi(z_1,z_2)=\frac{e^{i\theta}}{\sqrt{1+N}}z_1^{k_0}z_2^{N-k_0}.$$ 

It now remains to prove \eqref{eqnextremalD}. As $\Phi\in \mathcal D(\mathbb D^2)$ is inner, by Proposition~\ref{prop:RIF}, $\Phi$ is a rational inner function. Thus, there exists $F$ subset of $\mathbb T$ with $m(F)=0$ such that for $\zeta_1, \zeta_2 \in \mathbb T\backslash F,$ $\Phi(\zeta_1,\cdot)$ and $\Phi(\cdot,\zeta_2)$ are finite Blaschke products.
Note that for a finite Blaschke product $b$, by \cite[Proposition~3.5]{RS1991}, 
\begin{equation}\label{angular-derivative}
D_{\delta_\zeta}(b)=|b'(\zeta)|, ~~~\zeta\in \mathbb T.
\end{equation} 
Combining \eqref{eq:dmmphi} with \cite[Proposition~2.2]{RS1991} and \eqref{angular-derivative} , we have
\begin{align}
    D_{m,m}(\Phi)&=  \int_{\mathbb T\backslash F} \int_{\mathbb T}D_{\delta_\eta}(\Phi(\cdot, \zeta_2))dm(\eta)dm(\zeta_2) +\int_{\mathbb T\backslash F} \int_{\mathbb T}D_{\delta_\eta}(\Phi(\zeta_1,\cdot))dm(\eta)dm(\zeta_1)\notag \\
    &=  \int_{\mathbb T\backslash F} \int_{\mathbb T} |\partial_1\Phi(\eta,\zeta_2)|dm(\eta)dm(\zeta_2)+ \int_{\mathbb T\backslash F} \int_{\mathbb T} |\partial_2\Phi(\zeta_1,\eta)|dm(\eta)dm(\zeta_1)\notag\\
    &= \int_{\mathbb T\backslash F\times \mathbb T\backslash F} \left( |\partial_1\Phi(\zeta_1,\zeta_2)|+|\partial_2\Phi(\zeta_1,\zeta_2)| \right)dm(\zeta_1)dm(\zeta_2),
    \label{partial-1-int}
\end{align}
where $\partial_1\Phi(\zeta_1,\zeta_2)$ and  $\partial_2\Phi(\zeta_1,\zeta_2)$ denotes the angular derivatives of $\Phi(\cdot,\zeta_2)$ and $\Phi(\zeta_1,\cdot),$ respectively, at $\zeta_1,\zeta_2\in \mathbb T.$ Consider the function 
\begin{equation}\label{h-is-positive}
h(\zeta_1,\zeta_2)=|\partial_1 \Phi(\zeta_1,\zeta_2)| + |\partial_2\Phi(\zeta_1,\zeta_2)|,\quad \text{a.e.} (\zeta_1,\zeta_2)\in \mathbb T^2.
\end{equation}
By \eqref{partial-1-int}, $h\in L^1(\mathbb T^2).$

{\bf Claim:} $h$ is constant almost everywhere on $\mathbb T^2.$

{\bf Proof of Claim:}
 From \cite[Eqs.~4.14~\&~4.15]{M2013}, for an one-variable inner function $b$, we have
\begin{equation}\label{angular-d-bd-value}
\zeta  b'(\zeta)\overline{b(\zeta)} =|b'(\zeta)|\geqslant0,\quad \text{ a.e. } \zeta\in \mathbb T.
\end{equation}
As for each $\zeta_1,\zeta_2\in \mathbb T\backslash F,$ $\Phi(\cdot,\zeta_2)$ and $\Phi(\zeta_1,\cdot)$ are finite Blaschke products, $$h(\zeta_1,\zeta_2)=(\zeta_1 \partial_1 \Phi + \zeta_2 \partial_2 \Phi)(\zeta_1,\zeta_2)\overline{\Phi(\zeta_1,\zeta_2)},\quad \zeta_1,\zeta_2\in \mathbb T\backslash F.$$

As $\Phi$ is inner, so for each $\zeta_1,\zeta_2\in \mathbb T\backslash F,$  we know that $\Phi(\zeta_1,\cdot)$  and $\Phi(\cdot,\zeta_2)$ are finite Blaschke products. By Lemma~\ref{lem:d-mm-pq}, for any polynomials $p$ and $q$ we have
\begin{equation*}
\langle p\Phi,q\Phi\rangle_{\mathcal D(\mathbb D^2)}= \langle p,q \rangle_{\mathcal D(\mathbb D^2)}+ \int_{\mathbb T\backslash F\times \mathbb T\backslash F } p\overline{q}h(\zeta_1,\zeta_2)dm(\zeta_2)dm(\zeta_1).
\end{equation*}
Substituting $p=z_1^n$ and $q=z_1^kz_2^{\ell+1}$ for $n,k, \ell\geqslant 0,$  and using the fact that $\Phi$ is a constant multiple of an extremal function, we get
\begin{equation*}
0=\langle z_1^n\Phi,z_1^kz_2^{\ell+1}\Phi\rangle_{\mathcal D(\mathbb D^2)}= \int_{\mathbb T\backslash F\times \mathbb T\backslash F}\zeta_1^{n-k}(\overline{\zeta_2})^{\ell+1}h(\zeta_1,\zeta_2)dm(\zeta_2)dm(\zeta_1).
\end{equation*}
So $\hat{h}(a,b)=0$ for all $a\in \mathbb Z$ and $b\geqslant 1.$
Similarly substituting $p=z_2^n$ and $q=z_1^{k+1}z_2^{\ell}$,  for $n,k, \ell\geqslant 0,$ we obtain
 $\hat{h}(a,b)=0$ for all $b\in \mathbb Z$ and $a\geqslant 1.$ As $h$ is positive (real valued) (see \eqref{h-is-positive}) so $$\hat{h}(-a,-b)=\overline{\hat{h}(a,b)},\quad a,b\in \mathbb Z.$$ Thus we get $$\hat{h}(a,b)=0, \quad (a,b)\neq (0,0).$$ Therefore $h(\zeta_1,\zeta_2)$ is constant a.e. on $\mathbb T^2.$ Hence the claim.

Let say 
\begin{equation}\label{h=d}
h(\zeta_1,\zeta_2)=|\partial_1 \Phi(\zeta_1,\zeta_2)| + |\partial_2\Phi(\zeta_1,\zeta_2)|=N\quad \text{a.e. } (\zeta_1,\zeta_2)\in \mathbb T^2,
\end{equation}
 then 
$$|\partial_1 \Phi(\zeta_1,\zeta_2)|\leqslant N \quad  \text{and}\quad  |\partial_2\Phi(\zeta_1,\zeta_2)|\leqslant N,\quad \text{a.e. } (\zeta_1,\zeta_2)\in \mathbb T^2.$$

We know that for $\zeta_1,\zeta_2 \in \mathbb T\backslash F,$  $\Phi(\cdot,\zeta_2)$ and $\Phi(\zeta_1,\cdot)$ are finite Blaschke products, so they are holomorphic in a neighbourhood of $\overline{\mathbb D};$  $|\partial_1 \Phi(\zeta_1,\zeta_2)|\leqslant N,$ and $|\partial_2\Phi(\zeta_1,\zeta_2)|\leqslant N.$ Thus by Lemma~\ref{lem:partial-f-bounded}, $\partial_1\Phi, \partial_2\Phi\in H^\infty(\mathbb D^2)$.

For almost every $(\zeta_1,\zeta_2)\in \mathbb T^2,$
\begin{eqnarray*}
(z_1\partial_1 \Phi+z_2\partial_2 \Phi)^*(\zeta_1,\zeta_2)&=&\zeta_1 \partial_1 \Phi(\zeta_1,\zeta_2)+ \zeta_2 \partial_2 \Phi(\zeta_1,\zeta_2)\\
&\overset{\eqref{angular-d-bd-value}}{=}&\Phi(\zeta_1,\zeta_2)|\partial_1 \Phi(\zeta_1,\zeta_2)|+ \Phi(\zeta_1,\zeta_2)|\partial_2 \Phi(\zeta_1,\zeta_2)|\\
&=& \Phi(\zeta_1,\zeta_2)h(\zeta_1,\zeta_2)\\
&\overset{\eqref{h=d}}{=}& N \Phi(\zeta_1,\zeta_2).
\end{eqnarray*}
Since $\partial_1 \Phi,\partial_2 \Phi\in H^\infty(\mathbb D^2),$ the uniqueness of boundary functions yields \eqref{eqnextremalD}. This completes the proof.
\end{proof}

\subsection{Left-inverse commuting submodules of \texorpdfstring{$\mathcal D(\mathbb D^2)$}{}}
In this section we characterise left-inverse commuting submodules of $\mathcal D(\mathbb D^2).$ 
\begin{proof}[{\bf Proof of Corollary~\ref{inv-d-m-m}}] 
Applying Theorem~\ref{inv-sub-d-mu} we get that a submodule $\mathcal S$ of $\mathcal D(\mathbb D^2)$ is left inverse commuting if and only if $$\mathcal S=\varphi\mathcal D(\nu_1,\nu_2),$$ where $\varphi$ is an extremal function in $\mathcal D(\mathbb D^2)$ and $\nu_1,\nu_2$ has the form \eqref{nu-1-2-formula}. By Theorem~\ref{extremals-are-monomials}, extremal functions of $\mathcal D(\mathbb D^2)$ are constant multiples of monomials. So a submodule $\mathcal S$ of $\mathcal D(\mathbb D^2)$ is left inverse commuting if and only if $$\mathcal S=\varphi \mathcal D(\nu_1,\nu_2),$$ where $$\varphi(z_1,z_2)=\frac{e^{i\theta}}{\sqrt{1+k+\ell}}z_1^kz_2^\ell$$  for some $k,\ell\geqslant 0$ and $\theta\in \mathbb R,$ and 
$$d\nu_1(\zeta)=d\nu_2(\zeta)=\frac{1}{1+k+\ell} dm(\zeta).$$
Hence the result.
\end{proof}

\begin{remark}\label{non-monomial-example} 
Although all extremal functions in $\mathcal D(\mathbb D^2)$ are scalar multiple of monomials, this is not true for general $\mathcal D(\mu_1,\mu_2)$. For example, let $\lambda\in \mathbb D\backslash\{0\}$ and consider $\psi(w)=\frac{w-\lambda}{1-\overline{\lambda}w},$ $w\in \mathbb D.$ A straightforward computation yields that the function $\varphi(z_1,z_2):=\psi(z_1)/\sqrt{2}$ is extremal in $\mathcal D(|\psi'|^{-1}m,m)$.
\end{remark}

Guo and Zhao \cite[Theorem~1.1]{GZ2010} showed that two submodules of the classical Dirichlet space $\mathcal D(\mathbb D)$ are unitarily equivalent if and only if they are equal. 
In the following corollary, we characterize unitary equivalence of left-inverse commuting submodules in terms of the total degree of their monomial generators.

\begin{corollary}
Let $\mathcal S_1=[z_1^a z_2^b]$ and $\mathcal S_2=[z_1^c z_2^d]$ be left-inverse commuting submodules of $\mathcal D(\mathbb D^2)$, where $a,b,c,d$ are nonnegative integers. Then $\mathcal S_1$ and $\mathcal S_2$ are unitarily equivalent if and only if $a+b=c+d$.
\end{corollary}
\begin{proof}
Set $s=a+b+1$ and $t=c+d+1$. By Corollary~\ref{inv-d-m-m}, $\mathcal S_1$ and $\mathcal S_2$ have the representations
$$\mathcal S_1=\frac{z_1^a z_2^b}{\sqrt{s}}\,\mathcal D\left(\frac{m}{s},\frac{m}{s}\right)\text{ and }\mathcal S_2=\frac{z_1^c z_2^d}{\sqrt{t}}\, \mathcal D\left(\frac{m}{t},\frac{m}{t}\right).$$
By Theorem~\ref{thm:uni-equiv}, $\mathcal S_1$ and $\mathcal S_2$ are unitarily equivalent if and only if their model measures agree. This holds precisely when $s=t$, or equivalently, $a+b=c+d$.
\end{proof}

The following example shows that not all submodules of $\mathcal D(\mathbb D^2)$ are left-inverse commuting.
\begin{example}\label{not-LIC}
Let $\mathcal S=\{f\in \mathcal D(\mathbb D^2): f(0,0)=0\}.$ We claim that $z_1,z_2\in \mathcal S\ominus (z_1\mathcal S+z_2\mathcal S).$ Clearly, $z_1,z_2\in \mathcal S.$ Let $f,g\in \mathcal S.$ Then $f(0,0)=g(0,0)=0.$ Since monomials are orthogonal in $\mathcal D(\mathbb D^2)$,
\begin{align*}
\langle z_1,z_1f+z_2g\rangle_{\mathcal D(\mathbb D^2)}&=\langle z_1,z_1f\rangle_{\mathcal D(\mathbb D^2)}+\langle z_1, z_2g\rangle_{\mathcal D(\mathbb D^2)}\\
&= \big\langle z_1,\sum_{\substack{k,\ell\geqslant 0\\k+\ell\geqslant 1}}\hat{f}(k,\ell)z_1^{k+1}z_2^\ell\big\rangle_{\mathcal D(\mathbb D^2)}+\big\langle z_1, \sum_{\substack{k,\ell\geqslant 0\\k+\ell\geqslant 1}}\hat{g}(k,\ell)z_1^{k}z_2^{\ell+1}\big\rangle_{\mathcal D(\mathbb D^2)}\\
&= 0. 
\end{align*}
Similarly one can check that $\langle z_2,z_1f+z_2g\rangle_{\mathcal D(\mathbb D^2)}=0.$ Thus $z_1,z_2\in \mathcal S\ominus (z_1\mathcal S+z_2\mathcal S)$ and hence $\dim \mathcal S\ominus (z_1\mathcal S+z_2\mathcal S)\neq 1.$ So by Theorem~\ref{inv-sub-d-mu}(i), $\mathcal S$ is not left-inverse commuting.
\eof
\end{example}
We conclude this section with the following observation that there are no nonzero doubly commuting submodules of $\mathcal D(\mathbb D^2).$
\begin{theorem}\label{non-d-commuting}
There is no nonzero doubly commuting submodule of $\mathcal D(\mathbb D^2).$
\end{theorem}
\begin{proof}
Let $\mathcal S$ be a doubly commuting submodule of $\mathcal D(\mathbb D^2).$ Denote $T_i=M_{z_i}|_\mathcal S,$ $i=1,2.$ Then $(T_1,T_2)$ is a toral $2$-isometry, that is
\begin{equation*}
I - T_i^*T_i - T_j^*T_j+T_i^*T^*_jT_iT_j=0, \quad i,j = 1,2,
\end{equation*} 
\eqref{C1} 
and $T_1T_2^*=T_2^*T_1.$ Hence
\begin{align*}
0=I-T_1^*T_1-T_2^*T_2+T_1^*T_2^*T_1T_2&=I-T_1^*T_1-T_2^*T_2+T_1^*T_1T_2^*T_2=(I-T_1^*T_1)(I-T_2^*T_2).
\end{align*}
Since $T_i$ is $2$-isometry and hence expansive (see \cite[Lemma~1(a)]{R1988}) so $T_i^*T_i-I\geqslant 0$ for  $i=1,2.$ Again using the fact that $(T_1,T_2)$ is doubly commuting so 
\begin{equation*}
(T_1^*T_1-I)(T_2^*T_2-I)=(T_2^*T_2-I)(T_1^*T_1-I)=0.
\end{equation*}
Let $f\in \mathcal S\backslash\{0\}.$ Then for $i=1,2,$
\begin{equation}\label{defect-norm}
\langle (T_i^*T_i-I)f, f\rangle_{\mathcal D(\mathbb D^2)}=\langle z_if, z_if\rangle_{\mathcal D(\mathbb D^2)}-\langle f,f\rangle_{\mathcal D(\mathbb D^2)}= \|f\|^2_{H^2(\mathbb D^2)} >0.
\end{equation}
Hence $(T_i^*T_i-I)$ is injective for $i=1,2.$ Then $(T_2^*T_2-I)(T_1^*T_1-I)=0$ forces $T_1^*T_1-I=0.$
Which is a contradiction to \eqref{defect-norm}. Hence $\mathcal S=\{0\}.$
\end{proof}

\vskip .2cm
 \noindent {\bf Author Contributions} Both authors made equal contributions. 

\vskip .2cm
\noindent{\bf Data Availability} No datasets were generated or analysed during the current study.
\vskip .2cm
\noindent{\bf Competing Interests} The authors declare no competing interests.

\end{document}